\documentclass[11pt]{amsart}
\usepackage{amsmath}
\usepackage{amsxtra}
\usepackage{amscd}
\usepackage{amsthm}
\usepackage{amsfonts}
\usepackage{amssymb}
\usepackage{eucal}

\newtheorem{prop}{Proposition}

\newtheorem{conj}{Conjecture}
\newtheorem{thm}{Theorem}

\numberwithin{equation}{section}

\begin{document}

\newcommand{\thmref}[1]{Theorem~\ref{#1}}
\newcommand{\secref}[1]{\S~\ref{#1}}
\newcommand{\lemref}[1]{Lemma~\ref{#1}}
\newcommand{\propref}[1]{Proposition~\ref{#1}}
\newcommand{\corref}[1]{Corollary~\ref{#1}}
\newcommand{\conjref}[1]{Conjecture~\ref{#1}}
\newcommand{\remref}[1]{Remark~\ref{#1}}
\newcommand{\nc}{\newcommand}
\nc{\on}{\operatorname}
\nc{\Z}{{\mb Z}}
\nc{\C}{{\mb C}}
\nc{\F}{{\mb F}_q}
\nc{\Ql}{\bar{{\mb Q}}_l}
\nc{\bolda}{\mb A}
\nc{\cond}{|\,}
\nc{\bi}{\bibitem}
\nc{\pone}{{\mb P}^1}
\nc{\pa}{\partial}
\nc{\arr}{\rightarrow}
\nc{\al}{\alpha}
\nc{\ri}{\rangle}
\nc{\lef}{\langle}
\nc{\W}{{\mc W}}
\nc{\la}{\lambda}
\nc{\ep}{\epsilon}
\nc{\su}{\widehat{{\mf s}{\mf l}}_2}
\nc{\g}{{\mf g}}
\nc{\V}{{\mc V}}
\nc{\De}{\Delta}
\nc{\G}{\widehat{{\mf g}}}
\nc{\Li}{{\mc L}}
\nc{\out}{{\rm out}}
\nc{\ins}{{\rm in}}
\nc{\lo}{_{{\rm loc}}}
\nc{\crit}{{\rm crit}}
\nc{\Pro}{{\mb P}}
\nc{\z}{{\mc Z}}
\nc{\sw}{{\mf s}{\mf l}}
\nc{\wh}{\widehat}
\nc{\gb}{{\mf b}}
\nc{\n}{{\mf n}}
\nc{\h}{{\mf h}}
\nc{\kk}{h^\vee}
\nc{\zz}{{\mc Z}}
\nc{\wt}{\widetilde}
\nc{\hh}{\widehat{{\mf h}}}
\nc{\N}{\widehat{{\mf n}}}
\nc{\hi}{{\mc H}}
\nc{\M}{{\mc M}}
\nc{\Dh}{{\mc D}^{(1/2)}}
\nc{\larr}{\longrightarrow}
\nc{\La}{\Lambda}
\nc{\us}{\underset}

\nc{\mb}{\mathbb}
\nc{\mf}{\mathfrak}
\nc{\mc}{\mathcal}
\nc{\mbf}{\mathbf}

\title{Affine algebras, Langlands duality and Bethe Ansatz}

\author{Edward Frenkel}\thanks{\hspace*{-5mm} Published in the
Proceedings of the International Congress of Mathematical Physics,
Paris 1994, pp. 606-642. International Press, 1995\\ Partially
supported by NSF grant DMS-9205303.}

\address{Department of Mathematics, Harvard University, Cambridge,
MA 02138, USA}

\dedicatory{In memory of Claude Itzykson}

\maketitle

\setcounter{footnote}{0}

\section{Introduction.}
By Langlands duality one usually understands a correspondence between
automorphic representations of a reductive group $G$ over the ring of adels
of a field $F$, and homomorphisms from the Galois group
$\on{Gal}(\overline{F}/F)$ to the Langlands dual group $G^L$. It was
originally introduced in the case when $F$ is a number field or the field
of rational functions on a curve over a finite field \cite{La:fund}.

Recently A.~Beilinson and V.~Drinfeld \cite{BD:quant} proposed a
version of Langlands correspondence in the case when $F$ is the field
of rational functions on a curve $X$ over $\C$. This geometric
Langlands correspondence relates certain ${\mc D}$--modules on the
moduli stack $\M_G(X)$ of principal $G$--bundles on $X$, and
$G^L$--local systems on $X$ (i.e. homomorphisms $\pi_1(X) \arr
G^L$). A.~Beilinson and V.~Drinfeld construct this correspondence by
applying a localization functor to representations of the affine
Kac-Moody algebra $\G$ of critical level $k=-\kk$, where $h^\vee$ is
the dual Coxeter number.

The localization functor assigns a twisted ${\mc D}$--module on
$\M_G(X)$ to an arbitrary $\G$--module from a category ${\mc
O}^0$. The fibers of this ${\mc D}$--module are analogous to spaces
of conformal blocks from conformal field theory.  In fact, the ${\mc
D}$--module, which corresponds to the vacuum irreducible $\G$--module
of level $k \in \Z_+$, is the sheaf of sections of a vector bundle
(with projectively flat connection), whose fiber is dual to the space
of conformal blocks of the Wess-Zumino-Witten model. The analogy
between conformal field theory and the theory of automorphic
representations was underscored by E.~Witten in \cite{Witten}. It is
at the critical level where this analogy can be made even more precise
due to the richness of representation theory of $\G$.

The peculiarity of the critical level is that a completion of the
universal enveloping algebra of $\G$ at this level,
$U_{-\kk}(\G)=U(\G)/(K+\kk)$, contains a large center $Z(\G)$. This
center is isomorphic to the classical $\W$--algebra $\W(\g^L)$
associated to the Lie algebra $\g^L$, which is Langlands dual to $\g$
\cite{FF:ctr}. Recall that $\W(\g^L)$ consists of functionals on a
certain Poisson manifold, which is obtained by the Drinfeld-Sokolov
reduction from a hyperplane in the dual space to $\widehat{\g^L}$
\cite{DS}. Elements of this Poisson manifold can be considered as
connections of special kind on a $G^L$--bundle over a punctured disc
called $\g^L$--opers in \cite{BD:opers}. For example, $\sw_2$--opers
are the same as projective connections. Thus, $\g^L$--opers define
characters of $Z(\G)$, and one can obtain an infinite-dimensional
family of $\G$--modules of critical level by factoring Verma modules
by these characters. The $\G$--modules, which one obtains this way,
play the role of unramified representations of the group $G$ over a
local field, while $\g^L$--opers play the role of their local
Langlands parameters.

By applying localization functor to these $\G$--modules, one can
associate a ${\mc D}$--module on $\M_G(X)$ to each regular
$\g^L$--oper on $X$. On the other hand, a regular $\g^L$--oper defines
a $G^L$--local system on $X$ and hence a homomorphism $\pi_1(X) \arr
G^L$. This way one can establish, at least partially, the geometric
Langlands correspondence \cite{BD:quant}.

When $X$ is a rational curve, this correspondence can be analized
rather explicitly. In this case we should consider the moduli space of
$G$--bundles with parabolic structures at a finite number of marked
points. One can then attach a ${\mc D}$--module on this moduli space
to a $\g^L$--oper with regular singularities at the marked points. For
example, for $G=SL_2$ the ${\mc D}$--module thus obtained describes a
system of differential equations: $H_i \Psi = \mu_i \Psi$, where
$\mu_i$'s are the accessory parameters (residues at the marked points)
of the corresponding projective connection, and $H_i$'s are certain
mutually commuting differential operators. These operators can be
identified with the hamiltonians of the Gaudin model, which is a
completely integrable quantum spin chain associated to $\g$
\cite{G}. The connection between Gaudin's model and affine algebras
was studied in my joint work with B.~Feigin and N.~Reshetikhin
\cite{FFR}. There we gave a new interpretation of the Bethe ansatz.

Bethe ansatz is a method of diagonalization of commuting hamiltonians,
which is widely used in various models of statistical mechanics. This
method is one of the cornerstones of the powerful Quantum Inverse
Scattering Method, which was intensively developed in the last twenty
years, most notably, by the L.D.~Faddeev School, cf. e.g. \cite{Faddeev}.

The idea of Bethe ansatz is to look for eigenvectors in a particular
form. One can show that a vector written in such a form is an eigenvector
if a certain system of algebraic equations, called the Bethe ansatz
equations, is satisfied. The problem is to prove completeness of the Bethe
ansatz, which means that all eigenvectors can be written in this special
form.

For $G=SL_2$, I show in Sect. 5 below that the existence of an
eigenvector of the Gaudin operators with given eigenvalues implies
that the corresponding projective connection generates trivial
monodromy representation of the fundamental group of $\pone$ without
the marked points. Remarkably, that is precisely equivalent to the
Bethe ansatz equations. This observation enables us to prove
completeness of Bethe ansatz in the Gaudin model associated to
$SL_2$. Similar results can be obtained for other groups. For an
alternative approach, cf. \cite{RV}.

Let us now go back to the geometric Langlands correspondence. In
\cite{Dr} V.~Drinfeld gave a beautiful construction of Langlands
correspondence for the group $GL_2$ over the function field of a curve
over $\F$, cf. also \cite{Laumon}. This construction is intrinsically
geometric and can be carried out for a curve $X$ over $\C$. Thus, one
obtains two completely different realizations of Langlands
correspondence for the group $SL_2$: one via the localization functor
\cite{BD:quant} and the other via the construction of \cite{Dr}. An
interesting question is to establish their equivalence.

In Sect.~6, I essentially do this in the case when the curve $X$ is
rational (these result is joint with B.~Feigin). Namely, I show that
the equivalence of the two constructions amounts to a separation of
variables in the Gaudin system, which was introduced by E.~Sklyanin
\cite{Skl:sep} as an alternative to Bethe ansatz. It is possible that
Sklyanin's separation of variables can be generalized to establish the
equivalence of the two constructions of the geometric Langlands
correspondence for curves of higher genus and groups of higher rank.

It is interesting that the Separation of Variables, possibly the most
powerful method of solving quantum integrable models \cite{Sklyanin},
appears as an ingredient of the geometric Langlands
correspondence. Thus, Langlands philosophy unifies affine algebras,
${\mc D}$--modules on curves and on moduli spaces of bundles, and
integrable systems. We believe that this is a part of an even richer
picture, in which these relations are ``deformed'' in various
directions. For instance, Sklyanin \cite{Skl:quant,Sklyanin} has found
a separation of variables for the XXZ model, which is a
$q$--deformation of a Gaudin model. This suggests that there should
exist a ``quantum Langlands correspondence'' between certain systems
of $q$-difference equations, in which the role of affine algebra is
played by the corresponding quantum affine algebra. The latter also
has a large center at the critical level \cite{RS,S,DE,FR}. Moreover,
elements of the spectrum of the center of a quantum affine algebra can
be viewed as $q$--difference operators \cite{FR}. Hopefully, quantum
affine algebras can bridge the gap between affine algebras and groups
over a local non-archimedian field.

It is worth mentioning that besides geometric Langlands correspondence for
complex curves, there are other intriguing new examples of Langlands
duality \cite{places}. We consider understanding this new duality pattern
as one of the most challenging problems in contemporary mathematical
physics.

\smallskip
\noindent{\bf Acknowledgements.} This paper reflects results of my
joint works with B.~Feigin and N.~Reshetikhin. I would like to thank
A.~Beilinson and V.~Drinfeld for sharing with me their unpublished
results and conjectures, and for illuminating discussions. I am deeply
grateful to D.~Kazhdan for conveying to me his insights on the theory
of automorphic representations, and for his comments on the first
draft of the paper. In writing the last section, I benefited from
valuable discussions with D.~Kazhdan, I.~Mirkovic, K.~Vilonen, and
M.~Thaddeus.

Finally, I would like to thank the Organizing Committee for the opportunity
to give a lecture in the inspiring atmosphere of the International Congress
in Paris.

\section{Localization functor.}    \label{functor}

\subsection{Setup.}    \label{loc}
Let $G$ be a connected simply-connected simple Lie group over $\C$, $\g$ be
its Lie algebra, and $\G$ be the corresponding affine algebra -- the
extension of the (formal) loop algebra $L\g = \g((t))$ by
one-dimensional center $\C K$ \cite{Kac}.

Fix a smooth projective curve $X$ over $\C$, a point $p$ on $X$ and a
formal coordinate near this point. Let ${\mc P}$ be a principal
$G$--bundle over $X$ with a trivialization on the formal disc around
the point $p \in X$. Denote by $\g_{\mc P} = {\mc P} \times_G \g$
the vector bundle associated to the adjoint representation of $G$. Let
$\g_{\out}^{{\mc P}}$ be the Lie algebra of sections of the bundle
$\g_{\mc P}$ over $X \backslash p$. There is a natural embedding
$\g_{\out}^{{\mc P}} \arr L\g$ obtained by expanding sections in
Laurent power series at $p$. It can be lifted to an embedding
$\g_{\out}^{{\mc P}} \arr \G$.

Consider the category ${\mc O}^0$ of $\G$--modules, on which the Lie
subalgebra $\g_{\ins} = \g[[t]] \subset \G$ acts locally finitely.
The modular functor assigns to a $\G$--module $M$ of an arbitrary
level $k$ from the category ${\mc O}^0$, the space of coinvariants
$H(\g,X,{\mc P},M) = M/\g_{\out}^{{\mc P}} M$, cf. \cite{Segal}. The
dual space to $H(\g,X,{\mc P},M)$ is called the {\em space of
conformal blocks}. A conformal block is thus a linear functional $f$
on $M$, which satisfies the {\em Ward identities}: $f(g \cdot x) = 0,
\forall g \in \g_{\out}^{{\mc P}}, \forall x \in M$.

The space $H(\g,X,{\mc P},M)$ depends on ${\mc P}$, and we want to
study all of these spaces simultaneously. There is a standard way of
doing that, due to Beilinson and Bernstein \cite{BB:loc}. One
constructs a {\em localization functor}, which assigns to a
$\G$--module $M$ a twisted ${\mc D}$--module $\Delta(M)$ on the
moduli space $\M_G(X)$ of $G$--bundles on $X$, such that the fiber of
$\Delta(M)$ at ${\mc P}$ is isomorphic to $H(\g,X,{\mc P},M)$
(cf. e.g. \cite{Borel,KS} for an introduction to the theory of ${\mc
D}$--modules). Essentially, this just means realizing the Ward
identities as differential equations on $\M_G(X)$, so that the space
$H(\g,X,{\mc P},M)$ gets interpreted as the space of local solutions
of those equations near ${\mc P} \in \M_G(X)$.\footnote{compare this
with minimal models, where conformal blocks are solutions of
differential equations on the moduli spaces of curves}

To define this functor, observe that the moduli space of principal
$G$--bundles with a trivialization on the formal disc around $p$ is
isomorphic to the homogeneous space $\M_G = LG/G_{\out}$, where $LG =
G((t))$ is the Lie group of $L\g$, and $G_{\out}$ is the Lie group of
the Lie algebra $\g_{\out}$ corresponding to the trivial bundle.

For any integer $k$ one can define a line bundle $\xi^k$ on $\M_G$
together with a homomorphism from $\G$ to the Lie algebra of
infinitesimal symmetries $\xi^k$, such that the central element $K$
maps to the constant function $k$. This gives us a homomorphism from
the algebra $U_k(\G) = U(\G)/(K-k)$ to the algebra $D'_k$ of global
differential operators on $\xi^k$.

Let ${\mc D}'_k$ be the sheaf of differential operators on
$\xi^k$. In fact, such a sheaf can be defined for any $k\in\C$, so
that $\G$ maps to its global sections, and $K$ maps to $k$. For any
$\G$--module $M$ of level $k$ we can now define a left ${\mc
D}'_k$--module on $\M_G$ as $\Delta'(M) = {\mc D}'_k \otimes_{U(\G)}
M$. One can show as in \cite{BB:loc}, \S~3, that the fiber of
$\Delta'(M)$ at ${\mc P} \in \M_G$ is indeed isomorphic to
$H(\g,X,{\mc P},M)$.

\subsection{Moduli space.}    \label{stack}
It is well-known that the set of isomorphism classes of principal
$G$--bundles on $X$ is isomorphic to the double quotient $G_{\ins}
\backslash LG / G_{\out}$, where $G_{\ins} = G[[t]]$. However, this set is
not the set of $\C$--points of an algebraic variety; the structure of
algebraic variety can only be given to a subset of semi-stable
$G$--bundles. One can cure this problem by considering the moduli stack
${\mc M}_G(X)$ of $G$--bundles on $X$. The precise definition of the
algebraic stack ${\mc M}_G(X)$ can be found, e.g., in \cite{BL}, \S~3. This
way we do not throw out any bundles: the set of isomorphism classes of
$\C$--points of ${\mc M}_G(X)$ is by definition the set of isomorphism
classes of all $G$--bundles on $X$. Moreover, ${\mc M}_G(X)$ is the stack
theoretic quotient of the scheme (of infinite type) $\M_G=LG/G_{\out}$ by
the algebraic group $G_{\ins}$.

The line bundle $\xi^k$ on $\M_G$ is $G_{\ins}$--equivariant, and hence it
descends to a line bundle $\eta^k$ on ${\mc M}_G(X)$. Denote by ${\mc
D}_k$ the sheaf of differential operators on $\eta^k$. Note that ${\mc
D}_k$ can be defined for any $k \in \C$.

Since $M \in {\mc O}^0$, the action of the Lie subalgebra $\g_{\ins}$ on
the module $M$ is locally finite. Therefore we can integrate the action of
$\g_{\ins}$ on $M$ to an action of its Lie group $G_{\ins}$. The ${\mc
D}'_k$--module $\Delta'(M)$ then becomes $G_{\ins}$--equivariant. Hence it
is a pull-back of a ${\mc D}_k$-module $\Delta(M)$ on the moduli space
${\mc M}_G(X)$, cf. \cite{BS}, \S~4.5.

The spaces of coinvariants $H(\g,X,{\mc P},M)$, which are naturally
isomorphic for different trivializations of ${\mc P}$ at $p$, are
isomorphic to the fiber of $\Delta(M)$ at ${\mc P}$.

\smallskip
\noindent{\em Examples.} Let $V_k$ be the vacuum Verma module over $\G$ of
level $k$, i.e. $V_k = U_k(\G) \otimes_{U(\g_{\ins})} \C$. Then
$\Delta(V_k) = {\mc D}_k$. One obtains more interesting ${\mc
D}$--modules by taking non-trivial quotients of $V_k$.

For instance, if $k$ is a positive integer, then $V_k$ contains a
unique singular vector. Let $L_k$ be the irreducible quotient of
$V_k$. This is an integrable $\G$--module. The spaces $H(\g,X,{\mc
P},M)$ have the same dimension for different ${\mc P}$'s. Hence
$\Delta(L_k)$ is the sheaf of sections of a vector bundle over
$\M_G(X)$ with a projectively flat connection, whose fiber is
$H(\g,X,{\mc P},M)$.\qed

\smallskip
\noindent{\em Remark.} The localization functor first appeared in
\cite{BB:kl} and \cite{BK} in the following situation. Let $B$ be the Borel
subgroup of $G$ and $G/B$ be the flag manifold. To each integral weight
$\la$ one can associate a line bundle $\xi_\la$ on $G/B$. There is a
surjective homomorphism from $U(\g)$ to the algebra of differential
operators on $\xi_\la$. This allows to define the localization functor from
the category ${\mc O}$ of $\g$--modules to the category of ${\mc
D}_\la$--modules, where ${\mc D}_\la$ is the sheaf of differential
operators on $\xi_\la$. It assigns to a $\g$ module $M$, the ${\mc
D}_\la$--module $\Delta(M) = {\mc D}_\la \otimes_{U(\g)} M$. This functor
is exact if $\la$ is dominant \cite{BB:kl}. The adjoint functor assigns to
a ${\mc D}_\la$--module, the $\g$--module of its global sections. For
example, if $L_\la$ is the irreducible representation with an integral
dominant weight $\la$, then $\Delta(L_\la)$ is just the sheaf of sections
of $\xi_\la$. According to the Borel-Weil theorem, the space of global
sections of $\xi_\la$ is isomorphic to $L_\la$.

The space $\M_G$ can be thought of as an analogue of the flag manifold for
the loop group, and hence the localization of $\G$--modules described above
is an analogue of the Beilinson-Bernstein localization.\qed

We want to apply the localization functor $\Delta$ to $\G$--modules of
critical level $k=-h^\vee$. But first we review representation theory at
this level.

\section{Representations of critical level.}
\subsection{Null-vectors and central elements.}    \label{borel}
Let us call a vector $x \in V_{-\kk}$ a {\em null-vector}, if $\g_{\ins}
\cdot v = 0$. It follows from the Kac-Kazhdan formula for the determinant
of the Shapovalov form \cite{KK} that $V_k$ can contain null-vectors other
than the highest weight vector $v_k$ only if $k=-h^\vee$. From now on let
us denote $V_{-\kk}$ by $V$ and $v_{-\kk}$ by $v$. Each null-vector $x \in V$
defines an endomorphism $e_x$ of $V$ by the formula $e_x(P v) = P x,
\forall P \in U(\G)$, and each endomorphism $e$ of $V$ defines a null
vector as $e \cdot v$. Hence the space ${\mf z}(\G)$ of null-vectors in
$V$ is isomorphic to End$_{\hat{\g}}(V)$.

Define the completion $\wt{U}(\G)$ of $U(\G)$ as the inverse limit of
$U(\G)/U(\G) (\g \otimes t^n\C[t]), n>0$. The action of $\wt{U}(\G)$ is
well-defined on all modules from the category ${\mc O}^0$. Using the
structure of vertex operator algebra on $V$ \cite{FZhu,FF:ctr}, we can
attach to each vector $P \in V$ a power series $Y(P,z) = \sum_{m\in\Z} P_m$
$z^m$ (for precise definition, cf. \cite{BFLM}). The coefficients of
these power series are elements of $\wt{U}_{-\kk}(\G) = \wt{U}(\G)/(K+\kk)$,
and altogether they span a Lie subalgebra $U_{-\kk}(\G)\lo$ of
$\wt{U}_{-\kk}(\G)$, which is called the local completion of $U_{-\kk}(\G)$
\cite{FF:ctr}. For example, for each $A \in \g$, we have: $Y((A \otimes
t^{-1}) v,z) \equiv A(z) = \sum_{n\in\Z} (A \otimes t^n) z^{-n-1}$. This
shows that $\G \subset U_{-\kk}(\G)\lo$.

Let $Z(\G)$ be the center of $U_{-\kk}(\G)\lo$. Since $\G \subset
U_{-\kk}(\G)\lo$, elements of $Z(\G)$ are central in $\wt{U}_{-\kk}(\G)$. It
is easy to show that if $x \in {\mf z}(\G)$, then all coefficients of
$Y(x,z)$ lie in $Z(\G)$. Furthermore, all elements of $Z(\G)$ can be
obtained this way, and the natural map $Z(\G) \arr {\mf z}(\G)$ defined by
the formula $c \arr c \cdot v$ is surjective \cite{FF:ctr}. In particular,
we see that the algebra ${\mf z}(\G) \simeq \on{End}_{\hat{\g}}(V)$ is
commutative.

\smallskip
\noindent{\em Example.} Let $\{ J_a \}, a=1,\ldots,\dim \g$, be a basis of
$\g$ orthonormal with respect to an invariant bilinear form. The vector $S
= \frac{1}{2} \sum_a (J_a \otimes t^{-1})^2 v$ lies in ${\mf z}(\G)$. The
coefficients $S_n$ of the power series
\begin{equation}    \label{sug}
Y(S,z) \equiv \sum_{n
\in \Z} S_n z^{-n-2} = \frac{1}{2} \sum_a :J_a(z)^2:,
\end{equation}
where the columns stand for the normal ordering, lie in $Z(\G)$. They are
called the {\em Sugawara operators}.\qed

\subsection{Symbols of null-vectors.}    \label{inv}
Recall that the center of $U(\g)$ is a free polynomial algebra with
generators $C^i$ of orders $d_i+1, i=1,\ldots,l$, where $d_1,\ldots,d_l$
are the exponents of $\g$. In particular, $C^1 = \frac{1}{2} \sum_a
(J_a)^2$. It is natural to try to generalize the formula for $S$ by taking
other generators $C^i, i>1$, and replacing each $J_a$ by $J_a \otimes
t^{-1}$. Unfortunately, this approach does not produce elements of
${\mf z}(\G)$ (and hence $Z(\G)$) in general. Some constructions have been
given for the affine algebras of types $A_n^{(1)}, B_n^{(1)}$ and
$C_n^{(1)}$ \cite{Hayashi}, but explicit formulas for elements of
${\mf z}(\G)$ are not known in general.

The standard filtration on $U(\G)$ induces a natural filtration on the
module $V$, such that gr $V$ is isomorphic to
$S^*(L\g/\g_{\ins})$. Moreover, gr ${\mf z}(\G)$ is isomorphic to the
space of $\g_{\ins}$--invariants of $S^*(L\g/\g_{\ins})$. Similarly,
gr $U(\g) \simeq S^*(\g)$, and it is known that the center of $U(\g)$
is isomorphic to the $\g$--invariant part of $S^*(\g)$. We can apply
our naive approach to the symbols $\bar{C}^i$'s of the central
elements $C^i$'s of $U(\g)$, i.e. replace $J_a$ by $J_a \otimes
t^{-1}$ in $\bar{C}^i$ and consider it as an element $\bar{S}^i$ of
$S^*(L\g/\g_{\ins})$. It is easy to check that each $\bar{S}^i$ is
$\g_{\ins}$--invariant. One can therefore ask whether each $\bar{S}^i$
can be lifted to some $S^i \in {\mf z}(\G)$.

The results stated in the next section give us the affirmative answer
to this question for an arbitrary affine Lie algebra.

\subsection{The structure of the center.}
Recall that ${\mf z}(\G)$ is a commutative algebra isomorphic to
End$_{\hat{\g}}(V)$. The Lie algebra of vector fields on the punctured disc
acts on $\G$ according to the formulas: $[t^n \pa_t, A \otimes t^m] = m A
\otimes t^{n+m-1}$. This action induces an action of the Lie algebra of
regular vector fields on the disc on $V$ if we put $t^n \pa_t \cdot v = 0,
n\geq 0$. In particular, the vector field $-t \pa_t$ provides a
$\Z$--grading on $V$, such that $\deg A \otimes t^n = -n$ and $\deg v =
0$. One can show that ${\mf z}(\G) \subset V$ is preserved by this action.

\begin{thm}[\cite{FF:ctr}]    \label{end}
There exist $S^1,\ldots,S^l \in {\mf z}(\G)$, such that $\deg S^i =
d_i+1$, and ${\mf z}(\G) \simeq \C[\pa_t^n S^i]_{i=1,\ldots,l;n\geq
0}$. In particular, $S^1 = \frac{1}{2} \sum_a (J_a \otimes t^{-1})^2 v$.
\end{thm}

The theorem implies that the center $Z(\G)$ of $U_{-\kk}(\G)\lo$ is
generated in the appropriate sense by the Fourier coefficients of the power
series $S^i(z) \equiv Y(S^i,z), i=1,\ldots,l$, corresponding to $S^i$'s. We
will give a more precise description of $Z(\G)$ using the Poisson
structure, which is defined as follows.

For any $x_1, x_2 \in Z(\G) \subset \wt{U}_{-\kk}(\G)$, let $\wt{x}_1,
\wt{x}_2$ be their liftings to $\wt{U}(\G)$. Then $[\wt{x}_1,\wt{x}_2] =
(K+\kk) \wt{y} + (K+\kk)^2 (\ldots)$, for some $\wt{y} \in \wt{U}(\G)$. Put
$\{ x_1,x_2 \} = y \equiv \wt{y}$ mod $(K+\kk)$. One can check that $y \in
Z(\G)$, and that this operation is a Poisson bracket on
$Z(\G)$.\footnote{This was observed by V.~Drinfeld following T.~Hayashi's
work \cite{Hayashi}.} The Lie algebra of vector fields on the punctured
disc acts on $U_{-\kk}(\G)\lo$, and this action preserves $Z(\G)$. Moreover,
the induced action on $Z(\G)$ is hamiltonian: the vector fields $-t^{n+1}
\pa_t$ act as $\{ S_n,\cdot \}$, where $S_n$ are the Sugawara operators
\eqref{sug}.

Recall that for any simple Lie algebra $\g$ one can define its
Langlands dual $\g^L$ as the Lie algebra whose Cartan matrix is the
transpose of the Cartan matrix of $\g$. The following description of
the center $Z(\G)$ of $U_{-\kk}(\G)\lo$ was conjectured by Drinfeld.

\begin{thm}[\cite{FF:ctr}]    \label{ctr}
$Z(\G)$ is isomorphic to the classical $\W$--algebra $\W(\g^L)$
associated to $\g^L$, as a Poisson algebra. This isomorphism is equivariant
with respect to the action of the Lie algebra of vector fields on the
punctured disc.
\end{thm}

Classical $\W$--algebras are defined via the Drinfeld-Sokolov hamiltonian
reduction \cite{DS}. Let us recall this construction.

\subsection{Hamiltonian reduction.}    \label{dsred}
Let $C(\g)$ be the space of connections on the trivial $G$--bundle over the
punctured disc (more precisely, on Spec $\C((t))$). Such a connection
is a first order differential operator $\pa_t + A(t)$, where $A(t) \in \g
\otimes \Omega^1 = \g \otimes \C((t)) dt$.

Let $(\cdot,\cdot)$ be the standard invariant bilinear form on $\g$
normalized so that the square of the maximal root is equal to $2$
\cite{Kac}. Denote by $r^\vee$ the maximal number of edges connecting
two vertices of the Dynkin diagram of $\g$. Using the form $(\cdot,\cdot)$
we can identify the space $C(\g)$ with the hyperplane in $\G^*$, which
consists of linear functionals on $\G$ taking the value $1/r^\vee$ on
$K$. This hyperplane is equipped with a canonical Poisson structure, which
is the restriction of the Kirillov-Kostant structure on $\G^*$. The space
of local functionals on $C(\g)$ is then a Poisson algebra.

We have a natural projection $\imath: C(\g) \arr (L\n_+)^*$, which can be
interpreted as the moment map with respect to the coadjoint action of the
Lie group $LN_+$. We now perform the hamiltonian reduction with respect to
the one-point orbit $\chi \in (L\n_+)^*$, such that $\chi(e_\al \otimes
t^n) = - 1$, if $\al$ is a simple root and $n=-1$, and $\chi(e_\al \otimes
t^n) = 0$ otherwise. This means that we take the inverse image of $\chi$
and take its quotient by the action of $LN_+$. It turns out that this
action is free \cite{DS}, and the quotient is an infinite-dimensional
affine space ${\mc C}(\g)$ which is isomorphic (non-canonically) to
$\C((t))^l$.

\smallskip
\noindent{\em Remark.} Following Drinfeld, consider the inverse image
${\mc C}'(\g)$ of the set of all $\chi' \in (L\n_+)^*$, such that
$\chi'(e_\al \otimes t^n) \neq 0$ for some $n\in\Z$, if $\al$ is simple,
and $\chi'(e_\al \otimes t^n) = 0$ otherwise. Then the loop group $LB_+$ of
the Borel subgroup $B_+ \subset G$ acts on ${\mc C}'(\g)$, and the
quotient ${\mc C}'(\g)/LB_+$ is isomorphic to ${\mc C}(\g)$. After this
identification, we obtain a natural action of vector fields on the
punctured disc on ${\mc C}(\g)$, which is mentioned in \thmref{ctr}. One
can show that with respect to this action ${\mc C}(\g)$ is isomorphic to
the product of the space of projective connections (cf. below) and
$\oplus_{i=2}^l \Omega^{d_i+1}$.\qed
\smallskip

The classical $\W$--algebra $\W(\g)$ is by definition the space of local
functionals on ${\mc C}(\g)$.

\smallskip
\noindent{\em Example.} Here is a more explicit description of
$\W(\sw_n)$. In this case each $LN_+$--orbit in $\imath^{-1}(\chi)$
contains a unique element of the form
\begin{equation}    \label{special}
\pa_t -
\begin{pmatrix}
0 & q_1(t) & q_2(t) & \hdots & q_{n-2}(t) & q_{n-1}(t) \\
1 & 0 & 0 & \hdots & 0 & 0\\
0 & 1 & 0 & \hdots & 0 & 0\\
\hdotsfor{6} \\
0 & 0 & 0 & \hdots & 1 & 0
\end{pmatrix},
\end{equation}
which can be viewed as the $n$th order differential operator
\begin{equation}    \label{diffn}
\pa_t^n - q_1(t) \pa_t^{n-2} - \ldots - q_{n-1}(t).
\end{equation}

The isomorphism of \thmref{ctr} identifies the local functionals $q_i(m)$ $
=$ $\int q_i(t)$ $ t^{m+i} dt$ with central elements $S^i_m = \int S^i(z)
z^{m+i} dz$ $\in Z(\widehat{\sw}_n)$.\qed

A similar description of the spaces ${\mc C}(\g)$ for classical Lie
algebras was given in \cite{DS}.

\subsection{Opers.}    \label{opersdef}
Following \cite{BD:opers}, we will call elements of ${\mc C}(\g)$,
$\g$--{\em opers} on the punctured disc. Similarly, one can define the
space ${\mc C}_+(\g)$ of {\em regular} $\g$--opers on the disc. Beilinson
and Drinfeld have generalized this notion to the case of an arbitrary curve
$X$ \cite{BD:opers}. Let $G^{{\rm ad}}$ be the group of inner automorphisms
of $\g$. A $\g$--oper on $X$ is a $G^{{\rm ad}}$--bundle over $X$ with a
connection and an additional structure, such that locally it can be
considered as an element of ${\mc C}_+(\g)$.

For example, an $\sw_n$--oper on a curve $X$ is a rank $n$ (holomorphic)
vector bundle $E$ over $X$ defined up to tensoring with a line bundle,
equipped with a full flag of subbundles $0 = E_0 \subset E_1 \subset \ldots
\subset E_n = E$ and a connection $\nabla: E \arr E \otimes \Omega^1$, such
that $\nabla \cdot E_i \subset E_{i+1} \otimes \Omega^1$ and the induced
maps $E_i/E_{i-1} \arr E_{i+1}/E_i \otimes \Omega^1$ are isomorphisms for
$i=1,\ldots,n-1$. Thus, $E_i/E_{i-1} \simeq E_1 \otimes \Omega^{i-1}$, and
after tensoring with the appropriate line bundle we can set $E_i =
\Omega^{i-(n+1)/2}$. As a $PGL_n$--bundle, $E$ is uniquely defined by these
properties, in particular, it does not depend on the choice of
theta-characteristic.

Over the punctured disc, the bundle $E$ can be trivialized, and the
connection $\nabla$ can be brought to the form \eqref{special}. Hence we
obtain an equivalent definition: an $\sw_n$--oper is an $n$th order
differential operator acting from $\Omega^{(-n+1)/2}$ to
$\Omega^{(n+1)/2}$, such that locally it has the form \eqref{diffn},
i.e. its principal symbol is equal to $1$, and the sub-principal symbol
vanishes.

\subsection{$\G$--modules associated to opers.}    \label{gmod}
Let $M_{\chi,k}$ be the Verma module over $\G$. Recall that $M_{\chi,k} =
U_k(\G) \otimes_{U(\wt{\gb}_+)} \C_{\chi},$ where $\wt{\gb}_+ = (\gb_+
\otimes 1) \oplus (\g \otimes t \C[[t]])$ is the standard Borel subalgebra
of $\G$, and $\chi \in \h^*$. Denote by $v_{\chi,k}$ the highest weight
vector of $M_{\chi,k}$. Let $\rho \in {\mc C}(\g^L)$ be a $\g^L$--oper on
the punctured disc. By \thmref{ctr}, $\rho$ defines a central character,
i.e. a homomorphism $\wt{\rho}: Z(\G) \arr \C$. We define a $\G$--module
$M^\rho_\chi = M_{\chi,-\kk}/$ Ker $\wt{\rho}$.

\smallskip
\noindent{\em Example.} The space ${\mc C}(\sw_2)$ consists of projective
connections on the punctured disc, i.e. operators of the form $\pa_t^2 -
q(t)$ acting from $\Omega^{-1/2}$ to $\Omega^{3/2}$, where $q(t) =
\sum_{n\in\Z} q_n t^{-n-2}$. This space is an $\Omega^2$--torsor.

The module $M^{q(t)}_\chi$ is by definition the quotient of $M_{\chi,-2}$
by the submodule generated by the vectors $[S_m-q_m] v_{\chi,-2}, m \in
\Z$. But we know that $S_m v_{\chi,-2} = 0$ for $m>0$ and $S_0 v_{\chi,-2}
= \frac{1}{2} \sum_a [J_a \otimes 1]^2 v_{\chi,-2} = \frac{1}{4}
\chi(\chi+2) v_{\chi,-2}.$ Therefore $M^{q(t)}_{\chi,-2} \neq 0$ if and
only if $q_m=0, m>0$, and $q_0 = \frac{1}{4}\chi(\chi+2)$. This means that
the projective connection $\pa_t^2 - q(t)$ has regular singularity at the
origin, and $q(t) = \frac{1}{4} \chi(\chi+2)/t^2 + \ldots.$\qed

Likewise, when $\g=\sw_n$, the quotient of the Verma module $M_{\chi,-\kk}$
defined by the differential operator \eqref{diffn} is non-zero if and only
if this operator has regular singularity at the origin and the most
singular terms are equal to the values of $S^i_0$ on the highest weight
vector of $M_{\chi,-\kk}$. In general, we also have an infinite-dimensional
family of $\G$--modules of critical level parametrized by $\g^L$--opers
with regular singularities (in the sense of \cite{BD:opers}).

Another family of $\G$--modules can be obtained from the module $V$. One
can show that Spec End$_{\hat{\g}}(V) \simeq {\mc C}_+(\g^L)$, the space
of regular $\g^L$--opers on the formal disc. Therefore each $\rho \in {\mc
C}_+(\g^L)$ defines a character $\bar{\rho}$ of End$_{\hat{\g}}(V)$. We set
$V^\rho = V/$ Ker $\bar{\rho}$.

\subsection{Wakimoto modules.} There is another very important family
of $\G$--modules of critical level -- the Wakimoto modules.

To motivate the definition of Wakimoto modules, recall that in the study of
the center $Z(\g)$ of $U(\g)$ an important role is played by the
Harish-Chandra homomorphism $Z(\g) \arr \C[\h^*]$. The image of $z \in Z$
is a polynomial on $\h^*$, whose value at $\la \in \h^*$ is equal to the
value of $z$ on the Verma module $M_\la$.

In order to generalize this construction to the affine case with $L\h^* =
\h^* \otimes \Omega^1$ instead of $\h^*$, we need $\G$--modules, which
depend on elements of $L\h^*$. The Verma modules can not be used as they
depend only on elements of $\h^*$.

Roughly, we would like take a one-dimensional representation of $L\h$,
extend it trivially to a $L\n_-$ and take the induced $L\g$--module. But
the Lie algebra $L\n_+$ would act freely on the induced module. Hence the
action of the completion $\wt{U}(\G)$ on such a module is not well-defined,
and we can not use them to construct an analogue of the Harish-Chandra
homomorphism.

However, we can define a $\G$--module, which is partially induced and
partially coinduced \cite{FF:wak}, so that the Lie algebra $\n_+ \otimes
t^{-1} \C[t^{-1}]$ acts freely and the Lie algebra $\n_+ \otimes \C[t]$
acts cofreely. On these modules the action of $\wt{U}(\G)$
well-defined. But the definition of these modules requires regularization
(normal ordering) and this leads to a shift in the level -- it becomes
$-\kk$. These modules are called Wakimoto modules, cf. \cite{F:talk} for a
review.

\smallskip
\noindent{\em Example}. If $\g=\sw_2$, then the Wakimoto module of critical
level has the following realization \cite{Wak}. Consider the Heisenberg
algebra $\Gamma$ with generators $a_n, a^*_n, n\in\Z$, and relations
$[a_n,a^*_m] = \delta_{n,-m}.$ Let $M$ be the Fock representation of this
algebra generated by a vector ${\mbf v}$, such that $a_n {\mbf v} = 0, n\geq
0; a^*_n {\mbf v} = 0, n>0.$ Put $a(z) = \sum_{n\in\Z} a_n z^{-n-1}, a^*(z)
= \sum_{n\in\Z} a^*_n z^{-n}.$ Let $\{ e,h,f \}$ be the standard basis of
$\sw_2$. For any formal Laurent power series $\chi(z) = \sum_{n\in\Z}
\chi_n z^{-n-1}$, define an action of $\su$ on $M$ by the formulas:
\begin{eqnarray}
e(z) = a(z), \quad \quad h(z) = -2 :a(z) a^*(z): +
\chi(z),\nonumber\\
f(z) = - :a(z) a^*(z) a^*(z): - 2 \pa_z a^*(z) + \chi(z)
a^*(z).\nonumber
\end{eqnarray}
The corresponding $\su$--module is the Wakimoto module; we denote it by
$W_{\chi(z)}$.\qed

Because of normal ordering, if we perform a change of coordinates,
$\chi(t)$ will transform not as a one-form, but as a connection $\pa_t +
\frac{1}{2}\chi(t): \Omega^{-1/2} \arr \Omega^{1/2}$ over the punctured
disc. Similarly, the Wakimoto modules over $\G$ are parametrized by
connections on a certain $H^L$--bundle, where $H^L$ is the dual group to
the Cartan subgroup $H$ of $G$. This space is an $\h^* \otimes
\Omega^1$--torsor. If we choose a coordinate $t$ on the disc, we can
identify it with $\h^* \otimes \C((t))$.

Now let us return to the center $Z(\su)$. As we know, it is generated by
elements $S_n$ given by formula \eqref{sug}. To describe the action of
$Z(\G)$ on $W_{\chi(t)}$, it is therefore sufficient to describe the action
of the operators $S_n$. The operator $S_n$ acts on $W_{\chi(t)}$ by
multiplication by some $q_n\in \C$. Set $q(t) = \sum_{n\in\Z} q_n
z^{-n-2}.$ The projective connection corresponding to the action of
$Z(\su)$ on $W_{\chi(t)}$ is then $\pa_t^2 - q(t)$. One can show by direct
computation that the relation between this projective connection and the
connection $\pa_t + \frac{1}{2}\chi(t)$ has the form:
\begin{equation}    \label{mt}
\pa_t^2 - q(t) = (\pa_t - \frac{1}{2}\chi(t))(\pa_t +
\frac{1}{2}\chi(t)),
\end{equation}
which means that
\begin{equation}    \label{miura}
q(t) = \frac{1}{4} \chi(t)^2 - \frac{1}{2} \pa_t \chi(t).
\end{equation}

\subsection{Miura transformation.}
Formula \eqref{miura} is the transformation law from $\pa_t +
\frac{1}{2}\chi(t)$ to $\pa^2_t - q(t)$, which is called the {\em Miura
transformation}. An equivalent way to define it is by saying that the Miura
transformation of the connection
$$\pa_t -
\begin{pmatrix}
\chi(t)/2 & 0\\
1 & -\chi(t)/2
\end{pmatrix}$$ on the trivial rank two vector bundle on the punctured disc
is the unique gauge equivalent connection of the form
$$\pa_t -
\begin{pmatrix}
0 & q(t)\\
1 & 0
\end{pmatrix}.$$

One can also say this in a coordinate independent way, suitable for an
arbitrary curve $X$. Recall that an $\sw_2$--oper on $X$ is a rank 2 bundle
$E$ which is a non-trivial extension $0 \arr \Omega^{1/2} \arr E \arr
\Omega^{-1/2} \arr 0$, and a connection $\nabla$, which defines an
isomorphism between $\Omega^{1/2} \subset E$ and $(E/\Omega^{1/2}) \otimes
\Omega^1$. Suppose there is another subbundle of rank one, $E' \subset E$,
which is preserved by the connection $\nabla$, i.e. $\nabla \cdot E'
\subset E' \otimes \Omega^1$. The Miura transformation is by definition the
transformation from the local system $(E',\nabla')$ to the oper
$(E,\nabla)$.\footnote{The subbundle $E'$ does not exist if $X$ is a
compact curve of genus $g>1$. However, an invariant subbundle $E'$ may
exist on an affine curve, e.g. on a punctured disc.}

In general, the action of the center $Z(\G)$ on the Wakimoto modules is
given by the generalized Miura transformation introduced in
\cite{DS}.\footnote{this is a hamiltonian map relating the Poisson
structures of the generalized KdV and mKdV equations} This fact follows
directly from our construction of the isomorphism between the center
$Z(\G)$ and the ${\mc W}$--algebra ${\mc W}(\g^L)$ \cite{FF:ctr},
cf. also \cite{FFR}. Thus, the affine analogue of the Harish-Chandra
homomorphism is the Miura transformation.

In the case of $\sw_n$ the Miura transformation can be described as
follows. Suppose that we have an oper $(E,\nabla)$ with another full flag
of subbundles $E'_1 \subset E'_2 \subset \ldots E'_n$, which is preserved
by the connection, i.e. $\nabla \cdot E'_i \subset E'_i \otimes \Omega^1$
for all $i=1,\ldots,n$. Such a flag defines a connection $\nabla'$ on the
direct sum $E'=\oplus_{i=1}^n E'_i/E'_{i-1}$. The Miura transformation is
the transformation from the local system $(E',\nabla')$ to the oper
$(E,\nabla)$.

Locally, in the basis induced by the flag $E_\bullet$ the connection
$\nabla$ can be represented by formula \eqref{special}, whereas in the
basis induced by the flag $E'_\bullet$, $\nabla$ can be represented as
\begin{equation}    \label{cartan}
\pa_t -
\begin{pmatrix}
\chi_1(t) & 0 & \hdots & 0 & 0\\
1 & \chi_2(t) & \hdots & 0 & 0\\
\hdotsfor{5} \\
0 & 0 & \hdots & \chi_{n-1}(t) & 0\\
0 & 0 & \hdots & 1 & \chi_n(t)
\end{pmatrix}.
\end{equation}
Thus, locally the Miura transformation is just the transformation from the
connection \eqref{cartan} to the unique gauge equivalent connection of the
form \eqref{special}, cf. \cite{DS}. This amounts to the following
splitting of the differential operator \eqref{diffn} into the product of
the first order linear operators:
$$\pa_t^n - q_1(t) \pa_t^{n-2} - \ldots - q_{n-1}(t) = (\pa_t -
\chi_1(t)) \ldots (\pa_t - \chi_n(t)).$$

\subsection{$\G$--modules and connections.}    \label{gc}
As we have seen, $\G$--modules of critical level are parametrized by
geometric data, namely, local systems on the formal punctured
disc. Moreover, many properties of these local systems have representation
theoretic meaning. Here is an example, which we will need in \S~5.

Let us look at the Miura transformation \eqref{miura} as the Riccati
equation on $\chi(t)$ with fixed $q(t)$. Assume that the projective
connection $\pa_t^2 - q(t)$ has at most regular singularity at the origin,
i.e. that $q(t) = \sum_{n\leq 0} q_n t^{-n-2}$. In that case the Riccati
equation \eqref{miura} can be considered as a system of algebraic equations
on the Fourier coefficients of $\chi(t) = \sum_{m\leq 0} \chi_m t^{-m-1}$:
\begin{equation}    \label{rec}
\frac{1}{4} \sum_{i+j=n} \chi_i \chi_j + \frac{n+1}{2} \chi_n = q_n, \quad
\quad n\leq 0.
\end{equation}

Note that the first equation is $q_0 = \frac{1}{4} \chi_0(\chi_0+2),$ which
coincides with the Harish-Chandra homomorphism for $\sw_2$, i.e. the
formula expressing the value $q_0$ of the Casimir operator of $U(\sw_2)$ on
the Verma module with highest weight $\chi_0$. This is a general property
of the Miura transformation.

Simple analysis of the system \eqref{rec} shows that it has exactly two
solutions if $q_0 \neq m(m+2)/2$ for any $m \in \Z$. Each of the
corresponding two solutions $\chi(t)$ of \eqref{miura} gives rise to a
solution $\phi(t) = e^{-\frac{1}{2} \int^t \chi(s) ds}$ of the equation
\begin{equation}    \label{phi}
(\pa_t^2 - q(t)) \phi = 0.
\end{equation}
Thus, if $q_0 \neq m(m+2)/4$ for any $m \in \Z$, then there are two
linearly independent solutions of \eqref{phi}, and the monodromy matrix is
semi-simple. In this case there are two equivalent ways to attach to each
solution $\chi(t)$ an irreducible representation of $\su$ of critical
level. One is to take the quotient of the Verma module $M_{\chi_0,-2}$ by
the central character defined by $\pa_t^2-q(t)$, and the other is to take
the Wakimoto module $W_{\chi(t)}$.

Now suppose that $q_0 = m(m+2)/4$ for some $m \in \Z_+$. In this case there
exists a unique solution of \eqref{miura} the form $\chi(t) = -(m+2)/t +
\ldots$, but in general this is the only solution representable as a formal
Laurent power series, and other solutions contain logarithmic terms. This
means that the monodromy matrix for such $q(t)$ is the Jordan block with
the eigenvalue $\pm 1$.

If we set $\chi_0=m \in \Z_+$, we can uniquely determine $\chi_n, -m\leq
n\leq -1$, from equation \eqref{rec}. But then we obtain the equation
$\sum_{i,j<0;i+j=-m-1} \chi_i \chi_j = 2 q_{-m-1}$. It imposes a
non-trivial condition on the coefficients of $q(t)$, which can be written
as $P_m(q_{-1},\ldots,q_{-m-1})=0$, where $P_m$ is a polynomial of degree
$m+1$ if we put $\deg q_i = -i$. For example, $P_0=q_{-1}$, $P_1 =
2q_{-1}^2-q_{-2}$, etc. A solution of \eqref{miura} of the form $\chi(t) =
m/t + \ldots$ exists if and only if the equation $P_m=0$ is satisfied. In
this case the projective connection $\pa_t^2 - q(t)$ has monodromy $\pm 1$
around $0$.

Now consider the Verma module $M_{m,-2}, m \in \Z_+$. Let $V_{m,-2}$ be the
quotient of $M_{m,-2}$ by the submodule generated by the vector $f(0)^{m+1}
v_{m,-2}$; in particular, $V_{0,-2}$ is the vacuum module $V_{-2}$. For
$q(t) = m(m+2)/4t^2 + \ldots$, denote by $V_m^{q(t)}$ the quotient of
$V_{m,-2}$ by the central character corresponding to $\pa_t^2-q(t)$.

\begin{prop}    \label{monodromy}
The module $V_m^{q(t)}$ is non-zero if and only if the projective
connection $\pa_t^2 - q(t)$ has monodromy $\pm 1$.
\end{prop}

\noindent{\em Proof.} Suppose that $\pa_t^2 - q(t)$ has monodromy $\pm
1$. Then there exists $\chi(t) = m/t + \ldots$, which solves the Riccati
equation \eqref{miura}. Consider the Wakimoto module $W_{\chi(t)}$. One can
show that there exists a homomorphism $M_{m,-2} \arr W_{\chi(t)}$, which
maps $v_{m,-2}$ to ${\mbf v}$ and $f(0)^{m+1} v_{m,-2}$ to $0$. Hence we
obtain a non-trivial homomorphism $V_{m,-2} \arr W_{\chi(t)}$. But the
center $Z(\su)$ acts on $W_{\chi(t)}$ according to the central character
$q(t)$. Therefore this homomorphism factors through the homomorphism
$V_m^{q(t)} \arr W_{\chi(t)}$, and hence $V_m^{q(t)} \neq 0$. The ``only
if'' part can proved in a similar fashion.\qed

If the projective connection $q(t) = m(m+2)/4t^2 + \ldots$ has
monodromy $\pm 1$, then there is a one-parameter family of solutions of the
Riccati equation \eqref{miura} of the form $\chi(t) = m/t + \ldots$, and
one solution of the form $\chi(t) = -(m+2)/t + \ldots$. Thus, Wakimoto
modules with central character $q(t)$ are parametrized by points of
$\pone$.

\smallskip
\noindent{\em Remark.} Consider the singular vector of the Virasoro algebra
of type $(1,m+1)$. It is an element of the universal enveloping algebra of
the subalgebra of the Virasoro algebra generated by $L_{-1},L_{-2},\ldots$,
depending on the central charge $c$. If we write it in terms of $L_i/c,
i=-1,\ldots,-m-1$, then all coefficients are polynomial in $c^{-1}$,
cf. \cite{DIZ}. The classical limit of this vector is by definition the
polynomial in $q_i \equiv L_i/c$, which is obtained from this expression by
setting $c^{-1}=0$. Using quantum Drinfeld-Sokolov reduction it is easy to
show that this polynomial is equal to $P_m$ (compare with \cite{DIZ}).\qed

These results can be generalized for an arbitrary affine algebra $\G$.

\section{Geometric Langlands correspondence.}

\subsection{Localization.} Let $\rho_p \in {\mc C}_+(\g^L)$ be a
$\g^L$--oper on the formal disc around a point $p\in X$, and $V^{\rho_p}$
be the $\G$--module of critical level introduced in \secref{gmod}. We want
to apply the localization functor from \secref{functor} to the module
$V^{\rho_p}$. It should give us a ${\mc D}_{-\kk}$--module
$\Delta(V^{\rho_p})$ on the moduli space ${\mc M}_G(X)$. In fact, the
square of the line bundle $\eta^{-\kk}$ is isomorphic to the canonical
bundle $\omega$ on ${\mc M}_G(X)$. Hence ${\mc D}_{-\kk}$ is the sheaf of
differential operators acting on $\omega^{1/2}$ \cite{BD:quant}. Because of
that we denote ${\mc D}_{-\kk}$ by $\Dh$.

\smallskip
\noindent{\em Remark.} The critical level is $-\kk$, because we use the
invariant bilinear form on $\g$, with respect to which the square of the
maximal root of $\g$ is equal to $2$. If we were using the Killing form,
the critical level would be $-1/2$.\qed
\smallskip

Note that any regular $\g^L$--oper on $X$ can be restricted to the
formal disc Spec $\C[[t_p]]$.

\begin{thm}[A.~Beilinson, V.~Drinfeld \cite{BD:quant}]    \label{bd}
The ${\mc D}^{(1/2)}$--module $\Delta(V^{\rho_p})$ is non-zero if and only
if $\rho_p$ is the restriction of a regular $\g^L$--oper $\rho$ on $X$. In
that case $\Delta(V^{\rho_p})$ is determined by $\rho$ and does not depend
on the position of $p$.
\end{thm}

According to \cite{BD:quant}, the ${\mc D}^{(1/2)}$--module
$\Delta(V^{\rho_p})$ describes a system of differential equations on
$\M_G(X)$, which is a quantization of Hitchin's system \cite{Hitchin},
cf. also \cite{Faltings}. Let $\hi_G(X)$ be the moduli stack of Higgs pairs
$\hi_G(X)$, which is isomorphic to the cotangent stack of the moduli of
$G$--bundles. Hitchin has defined a map $\nu: \hi_G(X) \arr \oplus_{i=1}^l
H^0(X,\Omega^{d_i+1})$ and showed that it defines a completely integrable
system on $\hi_G(X)$.

By a quantization of the Hitchin system one should understand a
commuting set of differential operators acting on a line bundle over
$\M_G(X)$, whose symbols coincide with Hitchin's hamiltonians. Such
differential operators can be constructed using the center $Z(\G)$ in
the following way \cite{BD:quant}. The algebra $U_{-\kk}(\G)$ and its
local completion map to the ring of global differential operators
acting on the line bundle $\xi^{-\kk}$ over $\M_G$,
cf. \secref{loc}. The central elements from $Z(\G)$ define
differential operators on $\M_G$, which commute with the right action
of $G_{\ins}$. Hence one obtains a homomorphism $\jmath: Z(\G) \arr
D^{(1/2)}$, where $D^{(1/2)}$ is the algebra of global differential
operators acting on the line bundle $\omega^{1/2}$ over $\M_G(X) =
\M_G/G_{\ins}$.

Using the description of the center $Z(\G)$ provided by \thmref{ctr} and
the description of the ring of regular functions on $\hi_G(X)$ provided by
Hitchin's theorem, Beilinson and Drinfeld \cite{BD:quant} prove that if $G$
is simply-connected, the homomorphism $\jmath$ is surjective. In
particular, $D^{(1/2)}$ is commutative. Furthermore, Spec $D^{(1/2)}$ is
canonically isomorphic to the space of regular $\g^L$--opers on $X$, and
the symbol map gives an isomorphism gr $D^{(1/2)} \simeq
\C[\hi_G(X)]$. Therefore elements of $D^{(1/2)}$ are indeed quantizations
of the Hitchin hamiltonians.

Now each $\g^L$--oper $\rho$ on $X$ provides a character of $D^{(1/2)}$,
and hence a ${\mc D}^{(1/2)}$--module $\Dh \otimes_{D^{(1/2)}} \rho$. This
${\mc D}$--module is isomorphic to $\Delta(V^{\rho_p})$ defined above. It
describes the system of differential equations
\begin{equation}    \label{system}
\xi \cdot \Psi = \rho(\xi) \Psi, \quad \quad \xi \in D^{(1/2)}.
\end{equation}

Thus, the construction outlined above associates to an arbitrary
$\g^L$--oper on $X$, a system of differential equations (a twisted
${\mc D}$--module) on $\M_G(X)$. Beilinson and Drinfeld put this
construction in the context of Langlands correspondence.

\subsection{The local Langlands correspondence.}    \label{ll}
The Langlands correspondence is a correspondence (bijection) between
two different sets of objects associated to a field $F$ and a
reductive connected algebraic group $G$. Originally, it was formulated
when $F$ is either a number field, or $F=\F(X)$, the field of rational
functions on a smooth projective curve $X$ defined over a finite field
$\F$, as a far-reaching generalization of abelian class field theory
(which corresponds to the simplest case when $G$ is the multiplicative
group) \cite{La:fund}. We will now briefly discuss some (elementary)
aspects of the Langlands correspondence. This and the following
subsections are written for non-experts and aim at giving them a very
rough idea of a certain general picture which we will subsequently
want to compare to the picture over $\C$. For more information we
refer the reader to excellent reviews \cite{reviews}.

For simplicity we will restrict ourselves with the case of the
function field $F=\F(X)$ and a (split connected) simple algebraic
group $G$ defined over $\F$. For a closed point $x \in X$ let ${\mc
O}_x$ be the completion of the local ring of $x$, i.e. ${\mc O}_x
\simeq {\mb F}_{q_x}[[t]]$, and ${\mc K}_x$ be its field of fractions,
i.e. ${\mc K}_x \simeq {\mb F}_{q_x}((t))$, where $q_x = q^{\deg x}$.

Let $\Ql$ be the algebraic closure of the field $\Ql$ of $l$--adic
numbers, where $l$ does not divide $q$. A representation of the group
$G_x=G({\mc K}_x)$ in a $\Ql$--vector space $\pi_x$ is called smooth
if the stabilizer of any vector of $\pi_x$ is an open subgroup of
$G_x$. A smooth irreducible representation is called {\em unramified}
if there exists a non-zero vector $v_x \in \pi_x$ that is invariant
with respect to the subgroup $K_x = G({\mc O}_x)$; such a vector is
then unique up to multiplication by $\Ql^\times$.

A complete description of the set of equivalence classes of
irreducible unramified representations of $G_x$ is known. They are
parametrized by semi-simple conjugacy classes in the Langlands dual
group $G^L$ over $\Ql$. For the definition of the Langlands dual
group, cf. \cite{La:fund,reviews}.\footnote{a categorical definition,
in the spirit of geometric Langlands correspondence, is given in
\cite{Gin}} Here we will only say that the root and weight lattices of
$G$ and $G^L$ are interchanged with the coroot and coweight lattices,
respectively.

The correspondence between representations of $G_x$ and conjugacy classes
in $G^L(\Ql)$ can be described via the so-called {\em principal series}
representations of $G_x$. It is easiest to define such a representation
for $G=GL_n$. Any semi-simple conjugacy class in $G^L(\Ql)=GL_n(\Ql)$
contains a diagonal matrix $y = \on{diag}(y_1,\ldots,y_n), y_i \in \Ql$,
defined up to permutation. To such a class we can associate a character
of the upper Borel subgroup
$$\chi_y(b) = (q_x^{1-n} y_1)^{\nu_x(b_{11})} (q_x^{2-n}
y_2)^{\nu_x(b_{22})} \cdots y_n^{\nu_x(b_{nn})},$$ where $b_{ii}$'s are the
diagonal entries of $b$, and $\nu_x$ is the standard norm on ${\mc
K}_x$. This character then defines an induced representation of $GL_{n,x}$
in the space of locally constant functions $f(g)$ on $GL_{n,x}$, such that
$f(b x) = \chi_y(b) f(g)$ for all $b \in B, f \in GL_{n,x}$. This is a
principal series representation. It is known that it contains exactly one
irreducible unramified  component, which depends only on the
conjugacy class of $y$. Thus one obtains a correspondence between conjugacy
classes in $G^L(\Ql)$ and $G_x$--modules, which is called the local
Langlands correspondence.

\subsection{Global Langlands correspondence.}    \label{glob}
Now suppose that we are given a set $\{ y_x \}_{x\in X}$ of conjugacy
classes of $G^L(\Ql)$ for all points of $X$. Then by the local Langlands
correspondence we can associate to each of them a $G_x$--module, $\pi_x$.

Recall that the ring of adels ${{\bolda}}$ of $F$ is the restricted
product of all completions of $F$: ${{\bolda}} = \prod'_{x\in X} {\mc
K}_x.$ Restricted means that we only consider collections $(a_x)_{x\in
X}$ for which $a_x \in {\mc O}_x$ for all but finitely many $x$. Note
that $F$ ``diagonally'' embeds into ${{\bolda}}$. The group
$G({{\bolda}})$ is the restricted product of the groups $G_x, x \in
X$, and we can consider its representation in the restricted tensor
product $\otimes'_{x\in X} \pi_x$ of the modules $\pi_x$. The latter
space is spanned by elements of the form $\otimes_{x\in X} w_x$, such
that $w_x$ is the $K_x$--invariant vector $v_x$ for all but
finitely many $x$. It is clear that the action of $G({{\bolda}})$ on
$\otimes'_{x\in X} \pi_x$ is well-defined.

The group $G({{\bolda}})$ naturally acts on the space
$C(G(F)\backslash G({{\bolda}}))$ of locally constant functions on
the quotient $G(F)\backslash G({{\bolda}})$. An irreducible
representation of $G({{\bolda}})$ is called {\em automorphic} if it
appears in the decomposition of $C(G(F)\backslash
G({{\bolda}}))$. Let us write such a representation as the restricted
tensor product $\otimes_{x\in X} \pi_x$. Then it is called {\em
unramified} if each $\pi_x$ is unramified. In that case the tensor
product of $K_x$--invariant vectors $\otimes_{x\in X} v_x$ is
right-invariant with respect to the compact subgroup $K = \prod_{x\in
X} K_x$. Hence it defines a function on the double quotient $G(F)
\backslash G({\bolda})/K$, which is called the {\em automorphic
function} corresponding to $\otimes_{x\in X} \pi_x$. Note that $G(F)
\backslash G({{\bolda}})/K$ is the set of isomorphism classes of
$G$--bundles on $X$.

A natural question is to describe all unramified automorphic
representations of $G({{\bolda}})$. In view of the local Langlands
correspondence, we can pose this question in the following way: {\em for
what sets $\{ y_x \}_{x\in X}$ of conjugacy classes of $G^L(\Ql)$ is the
corresponding representation $\otimes_{x\in X} \pi_x$ automorphic}?

Global Langlands correspondence gives (at least conjecturally) an answer to
this question. It describes the ``consistency conditions'' on the set $\{
y_x \}$, which make the representation $\otimes_{x\in X} \pi_x$
automorphic.

One approach to characterize them is analytic, using $L$--{\em
functions}. To each local factor $\pi_x$ and a finite-dimensional
representation of $G^L$ one can associate a local $L$--function. The
product of these functions over all points of $X$ is a global
$L$--function of $\otimes_{x\in X} \pi_x$. It is known that if
$\otimes_{x\in X} \pi_x$ is automorphic, then global $L$--function has
analytic continuation and satisfies a functional equation (much as
Dirichlet's $L$--functions) \cite{La:fund,JL,GJ}. It turns out that
for $GL_2$, the converse is also true: if the global $L$-function
associated to the two-dimensional representation (and its twists by
all continuous characters of $F^\times \backslash {{\bolda}}^\times$)
has these nice properties, then the representation $\otimes'_{x\in X}
\pi_x$ is automorphic \cite{JL}. A similar characterization of
automorphic representations has also been obtained for other groups,
cf. \cite{reviews}.

Another approach, which is closer to us, is to relate automorphic
representations of $G(\bolda)$ to representations of the {\em Galois
group} $\on{Gal}(\overline{F}/F)$. If we restrict ourselves with
unramified representations of $G(\bolda)$, we should consider the
maximal unramified quotient of the Galois group, which is isomorphic
to the fundamental group $\pi_1(X)$ of $X$. Actually, to be precise,
we should consider the Weil group of $F$, which is the inverse image
of $\Z \subset \wh{\Z} \simeq \on{Gal}(\overline{\mb F}_q/\F)$ under
the homomorphism $\pi_1(X) \arr \on{Gal}(\overline{\mb F}_q/\F)$, but
we will ignore this subtlety.

To each point $x\in X$ one can associate canonically the {\em Frobenius}
conjugacy class $\on{Fr}_x$ of $\pi_1(X)$, cf. e.g. \cite{Milne},
p. 292. So any homomorphism $\sigma: \pi_1(X) \arr G^L(\Ql)$ defines a
collection $\{ \sigma(\on{Fr}_x) \}_{x\in X}$ of conjugacy classes in
$G^L$. Now we can formulate the global Langlands conjecture.

\begin{conj}    \label{global}
An irreducible unramified representation $\otimes'_{x\in X} \pi_x$ is
automorphic if and only if there exists a continuous homomorphism
$\sigma: \pi_1(X) \arr G^L(\Ql)$, such that each $\pi_x$ corresponds
to the conjugacy class $\sigma(\on{Fr}_x)$ in the sense of local
Langlands correspondence.
\end{conj}

In particular, this conjecture implies that there is a correspondence
between the set of isomorphism classes of the irreducible unramified
automorphic representations of $G({\bolda})$ and the set of
isomorphism classes of the homomorphisms $\pi_1(X) \arr
G^L(\Ql)$. However, in general this correspondence is not one-to-one,
because there may be more than one homomorphism $\pi_1(X) \arr
G^L(\Ql)$ with the given set of the conjugacy classes
$\sigma(\on{Fr}_x)$ (on the automorphic side this should be reflected
by the multiplicity of the corresponding representation in
$C(G(F)\backslash G({{\bolda}}))$ being greater than one). It is known
that this can not happen for $G=GL_n$, and so in this case we should
have a bijection. But for general $G$ the picture is considerably more
complicated.

\subsection{Local Langlands correspondence over $\C$.}
The complex analogues of the groups $G_x$ and $K_x$ for the field
$F=\C(X)$ are the loop group $LG$ and its subgroup $G_{\ins}$. Define
the category ${\mc O}^0_{{\rm crit}}$ of unramified $\G$--modules of
critical level, which consists of the modules, on which the action of
$\g_{\ins}$ is locally finite, and which contain
$\g_{\ins}$--invariant vector. On such modules, the action of the Lie
algebra $\g_{\ins}$ can be integrated to an action of the Lie group
$G_{\ins}$.

The analogue of a conjugacy class in the group $G^L$ is a regular
$\g^L$--oper on the formal disc. (Note that the monodromy of a
$G^L$--local system on a punctured disc gives rise to a conjugacy
class in $G^L$; but of course all regular opers correspond to the
class of the unit.)

The analogue of local Langlands correspondence will be the following:
{\em each regular $\g^L$--oper $\rho_x$ on the formal disc defines an
irreducible $\G$--module of critical level}. This is the $\G$--module
$V^{\rho_x}$ from \secref{gmod}. Indeed, this module is irreducible,
the action of the Lie algebra $\g_{\ins} = \g \otimes \C[[t]]$ on it
is locally finite, and it contains a unique $G_{\ins}$--invariant
vector -- the projection of the generating vector of $V$. Moreover,
one can show that the modules $V^{\rho_x}$ exhaust the irreducible
objects of the category ${\mc O}^0_{{\rm crit}}$.

It is quite obvious that the Wakimoto modules are the analogues of
representations of the principal series. For example, as we have seen
in \secref{gc}, for each regular $\sw_2$--oper on the disc, i.e. a
projective connection $\pa^2_t - q(t)$, where $q(t) = \sum_{n\leq -2}
q_n t^{-n-2}$, there is a family of solutions $\chi(t)$ of the
corresponding Riccati equation \eqref{miura}. These solutions are
analogues of the character $\chi$ in the definition of principal
series representation. The Wakimoto modules $W_{\chi(t)}$
corresponding to them contain a unique unramified component isomorphic
to $V^{q(t)}$. So it depends on $q(t)$, but not on $\chi(t)$, just
like an unramified representation depends on the conjugacy class of
$\chi$, but not on $\chi$.

\smallskip
\noindent{\em Remark.} One can generalize the local Langlands
correspondence for the group $GL_n$ over a local field to include those
representations, which are not unramified, but contain a fixed vector with
respect to the Iwahori subgroup. Such representations correspond to a pair
of semi-simple and unipotent conjugacy classes $y$ and $n$ in $G^L(\Ql)$,
such that $y n y^{-1} = n^{q_x}$ \cite{KL} (for a general group $G$ one
needs to introduce additional parameters).

The analogue of the Iwahori subgroup over $\C$ is the standard Borel
subalgebra $\wt{\gb}_+$ of $\G$, cf. \secref{borel}. Given a $\g^L$--oper
$\rho$ on the disc with regular singularity in the sense of
\cite{BD:opers}, we can construct a $\G$--module of critical level
$M^\rho_\chi$, which has a vector invariant with respect to
$\wt{\gb}_+$. Here $\chi$ should be such that the most singular part of
$\rho$ is given by the central character on the Verma module $M_\chi$ over
$\g$, cf. \secref{gmod}. But the module $M^\rho_\chi$ is not necessarily
irreducible, and to specify its irreducible component, one has to fix an
extra datum. For instance, in the case of $\g=\sw_2$, one has to fix in
addition to a projective connection $\pa^2_t - q(t)$, a solution of the
equation \eqref{phi} (or \eqref{miura}).\qed

\subsection{Global Langlands correspondence for curves over $\C$.}
\label{ho}
Now suppose that we are given a $\g^L$--oper $\rho_x$ for each point $x\in
X$. We can assign to $\rho_x$ the $\G$--module of critical level
$V^{\rho_x}$. Let $\g({{\bolda}}) = \prod'_{x\in X} \g((t_x))$ be
the restricted product of the loop algebras corresponding to the points of
the curve, and $\G({{\bolda}})$ be its one-dimensional central extension,
whose restriction to each factor coincides with the standard extension. The
restricted tensor product $\otimes'_{x\in X} V^{\rho_x}$ is naturally a
$\G({{\bolda}})$--module. Now we want to define an analogue of the
automorphicity property for this module.

Recall that over a finite field an automorphic $G({{\bolda}})$--module is the
module that can be realized as a subspace of the space of functions
on $G(F)\backslash G({{\bolda}})$. Hence such a module defines an automorphic
function on $G(F)\backslash G({{\bolda}})/K$. Over the complex field we can
not expect $\otimes'_{x\in X} V^{\rho_x}$ to be realized in the space of
functions on $G(F)\backslash G({{\bolda}})$. But we may hope to associate to
it a multivalued function on $G(F)\backslash G({{\bolda}})/K = \M_G(X)$, or,
even better, a holonomic system of differential equations with regular
singularities (which this function satisfies). The corresponding (twisted)
${\mc D}$--module on $\M_G(X)$ can be constructed by localization of the
$\G({{\bolda}})$--module $\otimes'_{x\in X} V^{\rho_x}$.

The localization functor can be constructed in the same way as in
\secref{functor}. The difference is that in \secref{functor} we realized
$\M_G(X)$ as $G_{\out} \backslash LG/G_{\ins}$, using one point $p \in
X$. Hence we assigned a twisted ${\mc D}$--module on $\M_G(X)$ to a single
$\G$--module attached to this point. Now we realize $\M_G(X)$ as
$G(F)\backslash G({{\bolda}})/K$ using all points of $X$. Hence to construct
a ${\mc D}$--module, we have to attach a $\G$--module $M_x$ of some level
$k$ from the category ${\mc O}^0$ to each point of $X$.

Then we can assign a twisted ${\mc D}$--module on $G(F)\backslash
G({\bolda})$ to the $\G({{\bolda}})$--module $\otimes_{x\in X} M_x$
in the same way as
in \secref{functor}. The action of $\g_{\ins}$ on each $M_x$ can be
integrated to an action of the group $G_{\ins}$. Therefore our ${\mc
D}$--module on $G(F)\backslash G({{\bolda}})$ is $K$--equivariant and
descends to a twisted ${\mc D}$--module on $\M_G(X)$, which we denote by
$\wt{\Delta}(\otimes_{x\in X} M_x)$. One can show that
$\wt{\Delta}(\otimes_{x\neq p} V_k \otimes M)$ (i.e. we attach a
$\G$--module $M$ of level $k$ to the point $p \in X$ and the vacuum module
$V_k$ to all other points) is isomorphic to $\Delta(M)$ defined in
\secref{functor}.

Let us specialize to $k=-h^\vee$. For a given set $\{ \rho_x \}_{x\in X}$
of local $\g^L$--opers, let us call the $\G({{\bolda}})$--module
$\otimes_{x\in X} V^{\rho_x}$ {\em weakly automorphic} if
$\wt{\Delta}(\otimes_{x\in X} V^{\rho_x}) \neq 0$.

We can now state a weak version of global Langlands correspondence over
$\C$: {\em the $\G({{\bolda}})$--module $\otimes_{x\in X} V^{\rho_x}$ is
weakly automorphic if and only if there exists a globally defined regular
$\g^L$--oper $\rho$ on $X$, such that for each $x\in X$, $\rho_x$ is the
restriction of $\rho$ to a small disc around $x$}. This statement is
analogous to \thmref{bd}. Moreover, if such $\rho$ exists, then
$\wt{\Delta}(\otimes_{x\in X} V^{\rho_x})$ is isomorphic to
$\Delta(V^{\rho_p})$ from \thmref{bd} for any $p \in X$. Let us denote this
${\mc D}^{(1/2)}$--module on $\M_G(X)$ attached to $\rho$ by
$\Delta_\rho$.

On the other hand, one can attach to a regular $\g^L$--oper on $X$ a
$G^L$--bundle with a connection, which is automatically integrable since
$\dim X = 1$, and hence a homomorphism $\pi_1(X) \arr G^L$. Therefore we
obtain a correspondence between the homomorphisms $\pi_1(X) \arr G^L$
corresponding to monodromies of $\g^L$--opers on $X$ and ${\mc
D}$--modules on ${\mc M}_G(X)$, cf. \cite{BD:quant}.

\subsection{Hecke operators.}
In the case of a finite field, the automorphic function on the double
coset $G(F)\backslash G({{\bolda}})/K$ corresponding to an automorphic
representation $\otimes_{x\in X} \pi_x$ is an eigenfunction of the
Hecke operators. Moreover, its eigenvalues are given by traces of
$\sigma(\on{Fr}_x)$ on the finite-dimensional representations of
$G^L(\Ql)$, cf. \cite{reviews}.

It is possible to define analogues of Hecke operators over $\C$ as
certain operations (correspondences) on ${\mc D}$--modules on ${\mc
M}_G(X)$, cf. \cite{Laumon,BD:quant}. Using these operations one can
strengthen the statement of global Langlands correspondence by saying
that $\Delta_\rho$ is an ``eigensheaf'' with respect to these Hecke
correspondences, with ``eigenvalues'' given in terms of the
$G^L$--local system defined by $\rho$ \cite{BD:quant}.

More generally, one expects that for {\em any} homomorphism $\sigma:
\pi_1(X) \arr G^L(\C)$ (equivalently, a $G^L$--local system on $X$, or
a flat $G^L$--bundle over $X$), there exists a unique (up to
isomorphism) ${\mc D}$--module ${\mc F}_\sigma$ on $\M_G(X)$, which is
automorphic with respect to $\sigma$ in the sense that it is an
``eigensheaf'' of Hecke correspondences with ``eigenvalues'' given in
terms of $\sigma$. This is a general form of geometric Langlands
correspondence over $\C$, which is due to Drinfeld.

If $\sigma$ can be realized by a regular $\g^L$--oper $\rho$ on $X$, then
$\Delta_\rho$ provides a candidate for ${\mc F}_\sigma$. Although not
every $G^L$--local system can be realized by a regular $\g^L$-- oper, it
can be realized by a $\g^L$--oper with regular singularities at a finite
number of points. Given such an oper, one can construct a candidate for
${\mc F}_\sigma$ by some modification of the localization functor defined
above \cite{BD:quant}. Unfortunately, $\sigma$ can be realized by different
$\g^L$--opers, and so one has to check that the resulting ${\mc
D}$--module depends only on $\sigma$.

\section{Gaudin model and Bethe ansatz.}

\subsection{Geometric Langlands correspondence in genus $0$.} In this
section we will consider ${\mc D}$--modules $\Delta_\rho$ when $X$ has
genus zero. In this case we will have to generalize the correspondence by
allowing ramifications at a finite number of marked points. So, on the
``Galois side'' we consider $G^L$--opers on $\pone$ with regular
singularities at the marked points. On the ``automorphic side'', we
consider twisted ${\mc D}$--modules on the moduli space of $G$--bundles
with parabolic structures at the marked points. Recall that a parabolic
structure on a $G$--bundle ${\mc P}$ at $p\in X$ is a choice of Borel
subgroup in the fiber of ${\mc P}$ over $p$.

For convenience, let us choose a global coordinate $t$ on $\pone$, and let
the marked points be $z_1,\ldots,z_N$, and $\infty$. We first consider the
case of $\g=\sw_2$. On the ``Galois side'', we have projective connections
on $\pone$ with regular singularities at $z_1,\ldots,z_N$, and
$\infty$. Such a connection has the form $\pa_t^2 - q(t)$ with
\begin{equation}    \label{qt}
q(t) = \sum_{i=1}^N \frac{c_i}{(t-z_i)^2} + \sum_{i=1}^N
\frac{\mu_i}{t-z_i},
\end{equation}
where $\sum_{i=1}^N \mu_i = 0$. On the ``automorphic side'' we consider
systems of differential equations on the moduli space of rank two bundles
over $\pone$ with trivial determinant and with parabolic structures, i.e. a
choice of a line in the fiber over each $z_i$ and $\infty$. The set of such
bundles is the double coset $G_{\out}\backslash \prod_i LG_i/\prod_i
\wt{B}_i$, where $LG_i$ is the copy of the loop group of $G=SL_2$
associated to the $i$th marked point, $\wt{B}_i$ is its Borel subgroup, and
$G_{\out}$ is the Lie group $G$ over the ring of rational functions on
$\pone$, which are regular outside $z_1,\ldots,z_N$ and $\infty$.

The open subset $\M^{(N)}$ corresponding to the trivial bundle ${\mc O}
\oplus {\mc O}$ is isomorphic to $\left( \pone \right)^{N+1}/SL_2$: we
choose lines in the two-dimensional fibers over $z_1,\ldots,z_N$ and
$\infty$ up to a diagonal action of $SL_2$.

\subsection{Construction of ${\mc D}$--module.}    \label{construction}
Now introduce $\la_i, i=1,\ldots,N$ and $\infty$, such that
$\la_i(\la_i+2)/4 = c_i, i=1,\ldots,N,$ and $\la_\infty(\la_\infty+2)/4 =
\sum_{i=1}^N (c_i + z_i \mu_i).$ Let $q_i(t-z_i), i=1,\ldots,N$, and
$q_\infty(t^{-1})$ be the expansions of $q(t)$ given by \eqref{qt} around
$z_i, i=1,\ldots,N$, and $\infty$, respectively. Let $M_{\la_i}^{q_i}$ be
the corresponding $\su$--modules of critical level introduced in
\secref{gmod}.

\hskip-.1cm For $\la_1,\ldots,\la_N,\la_\infty \in \Z$, we have a line bundle
$\boxtimes_{i=1,\ldots,N ;\infty} {\mc O}(\la_i)$ over $\left( \pone
\right)^{N+1}$.  This line bundle is equivariant with respect to the
diagonal action of $G$. Denote by ${\mc L}$ the corresponding line
bundle on $\M^{(N)}$, and let ${\mc D}_{\la_i}$ be the sheaf of
differential operators on ${\mc L}$. A generalization of the
localization functor from \secref{functor} applied to the tensor
product
$\otimes_{i=1,\ldots,N;\infty} M_{\la_i}^{q_i}$ gives us a ${\mc
D}_{\la_i}$--module on $\M^{(N)}$, which we denote by
$\Delta(\la_i,\mu_i)$.

Let $H_i, i=1,\ldots,N$, be the following elements of the algebra
$U(\sw_2)^{\otimes N}$:
\begin{equation}    \label{hi}
H_i = \sum_{j\neq i} \frac{e^{(i)} f^{(j)} + f^{(i)} e^{(j)} + \frac{1}{2}
h^{(i)} h^{(j)}}{z_i-z_j},
\end{equation}
where $a^{(i)}$ stands for the element $1 \otimes \ldots \otimes a \otimes
\ldots \otimes 1$ with $a$ in the $i$th place. The algebra $U(\sw_2)$ maps
to the algebra of global sections of the sheaf of differential operators on
${\mc O}(\la)$ (this algebra is well-defined for any $\la$). The elements
$H_i$ commute with the diagonal action of $SL_2$, and hence the
corresponding differential operators lie in the algebra of global sections
of ${\mc D}_{\la_i}$. One can also check directly that the operators $H_i$
commute with each other. The following Proposition follows from Proposition
1 of \cite{FFR}.

\begin{prop}    \label{deltalai}
The ${\mc D}_{\la_i}$--module $\Delta(\la_i,\mu_i)$ is isomorphic to the
quotient of ${\mc D}_{\la_i}$ by the left ideal generated by $H_i-\mu_i,
i=1,\ldots,N$.
\end{prop}
\noindent{\em Remark.} An analogue of \thmref{bd} is that the ${\mc
D}_{\la_i}$--module obtained by localization of $\otimes_i M_{\la_i}^{q_i}$
is equal to $0$, unless there exists $q(t)$ of the form \eqref{qt}, such
that $q_i(t-z_i), i=1,\ldots,N$, and $q_\infty(t^{-1})$ are the expansions
of $q(t)$ around $z_i, i=1,\ldots,N$, and $\infty$, respectively.

Note also that we have: $\sum_{i=1}^N H_i = 0$ and $\sum_{i=1}^N \mu_i =
0$, and that mutual commutativity of the operators $H_i$ follows from
\propref{deltalai}.\qed

\subsection{Gaudin model.}    \label{gaud}
The ${\mc D}_{\la_i}$--module $\Delta(\la_i,\mu_i)$ describes the system
of differential equations on $\M^{(N)}$:
\begin{equation}    \label{goden}
H_i \Psi = \mu_i \Psi.
\end{equation}
This system of equations appears in the {\em Gaudin model} of statistical
mechanics.

Gaudin's model \cite{G} is a completely integrable quantum spin chain. We
associate to the $i$th site of a one-dimensional lattice, a
finite-dimensional $\sw_2$--module $V_{\la_i}$ of highest weight $\la_i \in
\Z_+$ and a complex parameter $z_i, i=1,\ldots,N$. The space of states is
the tensor product $\otimes_{i=1}^N V_{\la_i}$. The hamiltonians are the
mutually commuting operators $H_i, i=1,\ldots,N$, given by \eqref{hi},
which act on the space of states. Gaudin's model can be considered as a
degeneration of the XXX or XXZ Heisenberg magnetic chains,
cf. \cite{Skl:sep}. In the latter the symmetry algebra is the Yangian of
$\sw_2$ or the affine quantum algebra $U_q(\su)$, respectively, while in
Gaudin's model, it is the affine algebra $\su$.

One of the main problems arising in the Gaudin model is to find the joint
spectrum of the operators $H_i$ on $\otimes_{i=1}^N V_{\la_i}$. If we
realize the operators $H_i$ as differential operators on $\M^{(N)}$, then
the diagonalization problem can be expressed as a system of differential
equations \eqref{goden} on $\M^{(N)}$. It is easy to write it down
explicitly on an open subset of $\M^{(N)}$.

Let $U_i$ be the big cell on the $i$th copy of $\pone$ and let $x_i$ be a
coordinate on $U_i$. The differential operators on $U_i$ corresponding to
the basis elements of the $i$th copy of $\sw_2$ are
\begin{equation}    \label{sl2}
e^{(i)} = -x_i^2 \frac{\pa}{\pa x_i} + \la_i x_i, \quad h^{(i)} = 2x_i
\frac{\pa}{\pa x_i} - \la_i, \quad f^{(i)} = \frac{\pa}{\pa x_i}.
\end{equation}
Note that $\M^{(N)} \simeq \left( \pone \right)^N/B$, where $B$ is the
Borel subgroup of $SL_2$. The system \eqref{goden} can be considered as a
system of $B$--invariant equations on $\left( \pone \right)^N$. The
restriction of this system to the product of $U_i$'s, which is an open
subset of $\left( \pone \right)^N$, reads:
$$\sum_{j\neq i} \frac{1}{z_i-z_j} \left[ - (x_i-x_j)^2 \frac{\pa^2}{\pa x_i
\pa x_j} + (x_i-x_j)\left( \la_i\frac{\pa}{\pa x_j} - \la_j\frac{\pa}{\pa
x_i} \right) + \frac{\la_i \la_j}{2} \right] \Psi$$
$$ = \mu_i \Psi.$$

The ${\mc D}$--module $\Delta(\la_i,\mu_i)$ is non-trivial for {\em
arbitrary} ``eigenvalues'' $\mu_i$'s, i.e. locally one can find solutions
of the system above for any set of $\mu_i$'s. But these solutions
generically have non-trivial monodromies around the diagonals.

In contrast, in Gaudin's model one is only interested in eigenvectors,
which lie in the space of states $\otimes_{i=1}^N V_{\la_i}$. Such an
eigenvector provides a {\em polynomial} solution of the system
\eqref{system}, which has no monodromies around the diagonals. This
solution gives rise to a surjective homomorphism $\Delta(\la_i,\mu_i) \arr
{\mc L}$.

\subsection{Bethe ansatz.}    \label{argument}
The diagonalization problem for the Gaudin model (as well as for the
majority of statistical models) is usually solved by the {\em Bethe ansatz}
method. This method consists of the following.

There is an obvious eigenvector in $\otimes_{i=1}^N V_{\la_i}$: the tensor
product $|0\ri$ of the highest weight vectors $v_{\la_i}$ of the
$V_{\la_i}$'s. One constructs other eigenvectors by acting on this vector
by certain elementary operators, depending on auxiliary parameters. For
$w\in \C, w\neq z_i$, and $a \in \{ e,h,f \}$ put
\begin{equation}    \label{ft}
a(w) = \sum_{i=1}^N \frac{a^{(i)}}{w-z_i}.
\end{equation}
Now introduce Bethe vectors
\begin{equation}    \label{bethevector}
|w_1,\ldots,w_m\ri = f(w_1) \ldots f(w_m) |0\ri \in \otimes_{i=1}^N
V_{\la_i}.
\end{equation}
Let us compute the action of Gaudin's hamiltonians on such a
vector. It is convenient to pass to a family of operators
\begin{equation}    \label{su}
S(t) = \sum_{i=1}^N \frac{\la_i(\la_i+2)/4}{(t-z_i)^2} + \sum_{i=1}^N
\frac{H_i}{t-z_i}, \quad \quad t\in \C\backslash \{ z_1,\ldots,z_N \}.
\end{equation}
It is clear that the diagonalization problems for $H_i$'s and $S(t), t\in
\C$, are equivalent.

Note that $\Psi$ is an eigenvector of the operators $H_i$ with the
eigenevalues $\mu_i$ if and only if $S(t) \Psi = q(t) \Psi$, where $q(t)$
is given by \eqref{qt}. Explicit computation shows that this formula holds
for $\Psi = |w_1,\ldots,w_m\ri$ if and only if
\begin{equation}    \label{be}
\sum_{i=1}^N \frac{\la_i}{w_j-z_i} - \sum_{s\neq j} \frac{2}{w_j-w_s} = 0,
\quad \quad j=1,\ldots,m.
\end{equation}
These equations are called {\em Bethe ansatz equations}.

If these equations are satisfied, then the Bethe vector \eqref{bethevector}
is an eigenvector, and also a highest weight vector of weight
$\la_{\infty}$ in $\otimes_{i=1}^N V_{\la_i}$. In this case the eigenvalue
$q(t)$ can be represented as
\begin{equation}    \label{eigenvalue}
q(t) = \frac{1}{4} \chi(t)^2 - \frac{1}{2} \pa_t \chi(t),
\end{equation}
where
\begin{equation}    \label{conn}
\chi(t) = \sum_{i=1}^N \frac{\la_i}{t-z_i} - \sum_{j=1}^m
\frac{2}{t-w_j}.
\end{equation}

The standard hypothesis is {\em completeness} of Bethe ansatz, which means
that the Bethe vectors exhaust all the eigenvectors, which are highest
weight vectors.

Now let us look at this hypothesis from the point of view of geometric
Langlands correspondence. Consider the ${\mc D}_{\la_i}$--module obtained
by localization of the tensor product of $\su$--modules
$\otimes_{i=1,\ldots,N;\infty} V_{\la_i}^{q_i}$
, where $q_i$'s were
defined in \secref{construction}, and $V_{\la}^q$ was defined in
\secref{gc}. The fibers of this ${\mc D}_{\la_i}$ over all points of
$\M^{(N)}$ are isomorphic to the space of coinvariants ${\mc
H}=
(\otimes_{i=1,\ldots,N;\infty}
V_{\la_i}^{q_i})/\g_{\out}$. But
$
(\otimes_{i=1,\ldots,N;\infty}
V_{\la_i,-2})/\g_{\out} \simeq
(
\otimes_{i=1,\ldots,N;\infty}
V_{\la_i})^G$,
cf. e.g. \cite{FFR}. Therefore ${\mc H}$ is isomorphic to the quotient of
\newline
$(\otimes_{i=1,\ldots,N;\infty}
V_{\la_i})^G$ by the action of the
operators $H_i-\mu_i, i=1,\ldots,N$.

If there exists an eigenvector of the Gaudin operators $H_i$'s with the
eigenvalues $\mu_i$'s in $(
\otimes_{i=1,\ldots,N;\infty}
V_{\la_i})^G$,
then the space ${\mc H}$ is non-zero. But then $V_{\la_i}^{q_i} \neq 0$
for all $i=1,\ldots,N; \infty$. By \propref{monodromy}, this implies the
following result.

\begin{prop}    \label{hold}
If there is an eigenvector of the Gaudin hamiltonians in $\otimes_{i=1}^N
V_{\la_i}$ with the eigenvalues $\mu_i, i=1,\ldots,N$, then all solutions
of the differential equation
\begin{equation}    \label{diffeq}
\left( \pa^2_t - \sum_{i=1}^N
\frac{\la_i(\la_i+2)/4}{(t-z_i)^2} - \sum_{i=1}^N \frac{\mu_i}{t-z_i}
\right) \phi = 0
\end{equation}
on $\pone$ have monodromies $\pm 1$ around $z_1,\ldots,z_N$ and $\infty$,
and hence the monodromy representation of the projective connection
$\pa_t^2 - q(t)$ defining the ${\mc D}_{\la_i}$--module
$\Delta(\la_i,\mu_i)$ is trivial.
\end{prop}

But this property of \eqref{diffeq} is equivalent to the Bethe ansatz
equations! Indeed, if the equation \eqref{diffeq} has monodromies $\pm 1$,
then there is a solution of the form
\begin{equation}    \label{solution}
\phi(t) = \prod_{i=1}^N (t-z_i)^{-\la_i/2} \prod_{j=1}^m (t-w_j),
\end{equation}
for some $w_1,\ldots,w_m \in \C$. But then formulas \eqref{eigenvalue} and
\eqref{conn} must hold. The left hand side of \eqref{eigenvalue} has no
poles at the points $w_j$'s. Therefore there should be no poles at $w_j$'s
in the right hand side. Straightforward computation shows that this
condition is equivalent to Bethe ansatz equations \eqref{be}. In other
words, Bethe ansatz equations mean that the Miura transformation
\eqref{eigenvalue} ``erases'' the extra poles of the connection
$\pa_t-\chi(t)/2$ given by \eqref{conn}.

On the other hand, if Bethe ansatz equations \eqref{be} are satisfied for
some $w_1,\ldots,w_m$, then formula \eqref{eigenvalue} holds. Hence
\eqref{solution} is a solution of the equation \eqref{diffeq}. A linearly
independent solution can be constructed as $\phi(t) \int^t \phi(w)^{-2}
dw$, and it is clear that both of them have monodromies $\pm 1$.

Thus, if there is an eigenvector of the Gaudin hamiltonians with the
eigenvalues $\mu_i, i=1,\ldots,N$, then there exist $w_1,\ldots,w_m$ such
that equations \eqref{eigenvalue} and \eqref{be} are satisfied. But then we
can construct this eigenvector explicitly as the Bethe vector
\eqref{bethevector}. The completeness of Bethe ansatz now boils down to the
linear independence of these vectors, which can be shown by elementary
methods for generic values of $z_i$'s.

We conclude that Bethe ansatz describes those ${\mc D}_{\la_i}$--modules
$\Delta(\la_i,\mu_i)$ on $\M^{(N)}$, which admit a surjective homomorphism
to an invertible sheaf. It tells us that the $\sw_2$--opers, to which they
correspond, generate trivial monodromy representations
\newline
$\pi_1(\pone\backslash \{z_1,\ldots,z_N,\infty \})$ $\arr PGL_2$.

\subsection{Generalization to other Lie algebras.} Starting from a
$\g^L$--oper on $\pone$ with regular singular points $z_1,\ldots,z_N$ and
$\infty$, one can define a twisted ${\mc D}$--module on the moduli space
of $G$--bundles on $\pone$ with parabolic structures at $z_1,\ldots,z_N$
and $\infty$.

The open subset of the moduli space corresponding to the trivial
$G$--bundle is isomorphic to $F^{N+1}/G$. The restriction of our twisted
${\mc D}$--module to this subset is the quotient of the appropriate sheaf
of differential operators on $F^{N+1}/G$ by the left ideal generated by the
images of the central elements from $Z(\G)$. It is easy to write down a
formula for the generators of the ideal corresponding to the Sugawara
elements: we just have to replace the numerator in \eqref{hi} by $\sum_a
J_a^{(i)} J_a^{(j)}$. These operators are the hamiltonians of the Gaudin
model associated to $\g$. The new fact is that in general the subalgebra of
$U(\g)^{\otimes (N+1)}$ of elements commuting with $H_i$'s and among
themselves is quite large, because it includes elements corresponding to
higher order central elements from $Z(\G)$.

The existence of an eigenvector of the generalized Gaudin hamiltonians is
equivalent to the existence of a surjective homomorphism from the
corresponding ${\mc D}$--module to an invertible sheaf on $F^{N+1}/G$. A
straightforward generalization of our argument from \secref{argument}
implies that such a homomorphism exists only if the monodromy
representation of the corresponding $\g^L$--oper on $\pone$ is trivial. The
latter condition can be expressed as a system of algebraic equations
similar to the Bethe ansatz equations for $\g=\sw_2$. If these equations
are satisfied, one can construct explicitly the corresponding eigenvectors,
which look similar to the Bethe vectors for $\g=\sw_2$. The completeness of
Bethe ansatz then follows from linear independence of these vectors.

In \cite{FFR} eigenvectors of generalized Gaudin hamiltonians were
constructed explicitly using Wakimoto modules of critical level. Another
construction was given in \cite{RV}. The formula for these eigenvectors was
suggested in \cite{BF}, though the ``off-shell Bethe ansatz equations'' of
\cite{BF} had already been proved in \cite{Ch}.

Let us explain in more detail the case of $\g=\sw_3$; the general case is
similar.

\subsection{The case of $\sw_3$.} We start with a $\sw_3$--oper $\rho$,
which is a third order differential operator on $\pone$ with regular
singularities at $z_1,\ldots,z_N$, and $\infty$:
\begin{equation}    \label{oper}
\rho = \pa_t^3 - \sum_{i=1}^N \left(
\frac{c^1(\la_i)}{(t-z_i)^2} + \frac{\mu_i}{t-z_i} \right) \pa_t
- \sum_{i=1}^N \left( \frac{c^2(\la_i)}{(t-z_i)^3} +
\frac{\nu_i}{(t-z_i)^2} + \frac{\kappa_i}{t-z_i} \right).
\end{equation}

To construct the corresponding ${\mc D}$--module, we have to choose the
weights $\la_1,\ldots$, $\la_N$ and $\la_\infty$ of $\sw_3$, such that
$c^1(\la)$ and $c^2(\la)$ are the values of the central elements $C^1$ and
$C^2$ of $U(\sw_3)$ of orders $2$ and $3$, respectively, on the Verma
module $M_\la$. We will assume that $\la_1,\ldots,\la_N,\la_\infty$ are
integral dominant weights. Let
$$\rho_i = \pa^3_t - q^1_i(t-z_i) \pa_t + q^2_i(t-z_i), \quad \quad
i=1,\ldots,N,$$ be the expansions of the oper \eqref{oper} at $z_i$, and
$\rho_\infty$ be its expansion at $\infty$. Each $\rho_i$ defines a
character of the center $Z(\wh{\sw_3})$. Recall that in this case the
center is generated by the Fourier components of the Sugawara field
$S^1(z)$ of order $2$ and a field $S^2(z)$ of order $3$. The values of the
central character on them are given by the Fourier components of $q^1_i$
and $q^2_i$, respectively. Note that $c^1(\la_i)$ and $\mu_i$ are the
values on the components $S^1_0$ and $S^1_{-1}$, respectively, while
$c^2(\la_i), \nu_i$ and $\kappa_i$ are the values on the components $S^2_0,
S^2_{-1}$ and $S^2_{-2}$, respectively.

To each marked point we can now associate a $\G$--module of critical level,
$M_{\la_i}^{\rho_i}$, and apply the localization functor to the tensor
product of these modules. This gives us a ${\mc D}$--module on
$F^{N+1}/SL_3$. The corresponding system of differential equations consists
of $N$ second order differential equations, and $2N$ third order
differential equations.

The existence of an eigenvector of the generalized Gaudin hamiltonians in
$\otimes_{i=1}^N V_{\la_i}$ implies, in the same way as in the case
$\g=\sw_2$, that the monodromy representation $\pi_1(\pone\backslash \{
z_1,\ldots,z_N,\infty \}) \arr PGL_3$ defined by the oper \eqref{oper} is
trivial.

But the monodromy is trivial if and only if this oper admits a global Miura
transformation:
\begin{equation}    \label{dec}
\rho = (\pa_t - \chi_1(t))(\pa_t - \chi_2(t))(\pa_t - \chi_3(t)),
\end{equation}
where $\chi_i(t)$ are the components of a diagonal connection
$\pa_t+\chi(t)$.

To see that, observe that it follows from formulas \eqref{oper} and
\eqref{dec} that $\chi(t)$ must be of the form
\begin{equation}    \label{diagonal}
\chi(t) = \sum_{i=1}^N \frac{\la_i}{t-z_i} + \sum_{s\in W; s\neq 1}
\sum_{j_s=1}^{m_s} \frac{s(\bar{\rho})-\bar{\rho}}{t-w^{(s)}_{j_s}},
\end{equation}
where $W$ is the Weyl group of $\sw_3$. Here the weights of $\sw_3$ are
represented by diagonal matrices, so that the fundamental weights are
represented as $\omega_1 = \on{diag}[2/3,-1/3,-1/3]$ and $\omega_2 =
\on{diag}[1/3,1/3,-2/3]$.

Indeed, it is clear from \eqref{dec} that if $\chi(t)$ has a pole different
from $z_1,\ldots,z_N$, then it can not be of order greater than $1$. If
such an extra singular term of the form $\gamma/(t-w)$ exists for $w\neq
z$, then $c^1(\gamma)=c^2(\gamma)=0$, because $\rho$ is assumed to be
regular at $w$. But the weights which satisfy this condition lie in the
orbit of $0$ with respect to the action of the Weyl group $W$ of $\sw_3$
shifted by the half-sum of the positive roots $\bar{\rho}$. For $\g=\sw_3$
these weights are:
$0,-\al_1,-\al_2,-2\al_1-\al_2,-\al_1-2\al_2,-2\al_1-2\al_2$.

If $\rho$ can be represented in the form \eqref{dec} with $\la_i$'s
integral dominant weights, then it is easy to find linearly independent
solutions of the equation $\rho \cdot \phi = 0$: $$e^{\int^t \chi_3(u)
du},\quad \quad e^{\int^t \chi_3 du} \int^t dx e^{\int^x (\chi_2-\chi_3)
du},
$$
$$
e^{\int^t
\chi_3 du} \int^t dx \, e^{\int^x (\chi_2-\chi_3) du} \int^x dy \,
e^{\int^y (\chi_1-\chi_2) du}.$$ It follows from these formulas that in
the neighborhood of each point $z_i$ local solutions look as follows:
$(t-z_i)^{(-n^1_i-2n^2_i)/3}(1+\ldots), \quad
(t-z_i)^{(n^2_i-n^1_i)/3+1}(1+\ldots)$, and
$(t-z_i)^{(2n^1_i+n^2_i)/3+2}(1+\ldots),$ where we put $\la_i=n^1_i
\omega_1 + n^2_i \omega_2$. Hence the monodromy of each of these solutions
around $z_i$ is a third root of unity, i.e. the identity in $PGL_3$.

On the other hand, if $\rho$ has a trivial monodromy representation, then
it has three linearly independent global solutions, and hence it is
representable in the form \eqref{dec}.

\smallskip
\noindent{\em Remark.}  Recall that $\rho$ admits a Miura transformation if
and only if it admits a flag of invariant subbundles, i.e. the underlying
rank $3$ bundle $E$ has invariant subbundles $E'_1$ and $E'_2$ of ranks $1$
and $2$. But the oper structure already defines subbundles, $E_1$ and $E_2$
of ranks $1$ and $2$, cf. \secref{opersdef}. Thus, the fiber of $E$ at each
point of $\pone\backslash \{ z_1,\ldots,z_N,\infty \}$ contains two
flags. The function $\chi(t)$ has a singularity at a point $w$ if these two
flags are not in generic position at $w$. Non-generic positions of two
flags are labeled by elements of the Weyl group. Hence we obtain a
summation over the Weyl group in \eqref{diagonal}.\qed

\subsection{Generalized Bethe equations.} Representability of $\rho$ in the
form \eqref{dec} is equivalent to a system of algebraic equations on the
positions of the points $w^{(s)}_{j_s}$. Indeed, we have already made the
most singular terms of the right hand side of \eqref{dec} vanish, but
vanishing of other singular terms imposes additional relations. Suppose for
simplicity that in \eqref{diagonal} only the simple reflections $s_i$
occur. Since $s_i(\bar{\rho})-\bar{\rho} = -\al_i$, we have in this case:
\begin{equation}    \label{chi}
\chi(t) = \sum_{i=1}^N \frac{\la_i}{t-z_i} - \sum_{j=1}^m
\frac{\al_{i_j}}{t-w_j}.
\end{equation}
According to Proposition 6 in \cite{FFR} the corresponding equations read:
\begin{equation}    \label{beq}
\sum_{i=1}^N \frac{(\la_i,\al_{i_j})}{w_j-z_i} - \sum_{l\neq j}
\frac{(\al_{i_l},\al_{i_j})}{w_j-w_l} = 0, \quad j=1,\ldots,m.
\end{equation}
They can be called the generalized Bethe ansatz equations for
$\g=\sw_3$. Analogous equations can be written for an arbitrary Lie algebra
$\g$, cf. \cite{FFR}.

The equations \eqref{beq} have been interpreted in \cite{FFR} as the
equations on the existence of a null-vector (cf. \secref{borel}) in the
Wakimoto module $W_{\chi^j(t)}$, where $\chi^j(t-w_j)$ is the expansion of
$\chi(t)$ at $w_j$ for $i=1,\ldots,m$. If equations \eqref{beq} hold, then
such null-vectors exist and using them one can construct explicitly an
eigenvector of the Gaudin hamiltonians, cf. Theorem 3 in \cite{FFR}. The
formula for this eigenvector and the equations \eqref{beq} can also be
derived from the asymptotic analysis of solutions of the
Knizhnik-Zamolodchikov equation \cite{RV}, cf. also \cite{BF,Ch}.

\subsection{Remarks.}
1. Formula \eqref{solution} shows that in the case $\g=\sw_2$
the numbers $w_j$'s defining the Bethe vectors can be found as the zeroes
of the polynomial $\phi \prod_{i=1}^N (t-z_i)^{\la_i/2}$, where $\phi$ is a
solution of the differential equation \eqref{diffeq} (the one for which
this expression is really a polynomial).

In general one can find the numbers $w_j$'s in a similar way. Let us
explain this for $\g=\sw_3$ in the case when we only have degeneracies
corresponding to the simple reflections. Then the positions of the
degeneracies corresponding to the simple reflection $s_2$ are the zeroes of
the polynomial $\phi \prod_{i=1}^N (t-z_i)^{(n^1_i+2n^2_i)/3}$, where
$\phi$ is a solution of the equation $\rho \cdot \phi = 0$. On the other
hand, the positions of the degeneracies corresponding to the simple
reflection $s_1$ are the zeroes of the polynomial $\wt{\phi} \prod_{i=1}^N
(t-z_i)^{(2n^1_i+n^2_i)/3}$, where $\wt{\phi}$ is a solution of the
equation $\rho^{{\rm ad}} \cdot \wt{\phi} = 0$. Here $\rho^{{\rm ad}}$
stands for the adjoint of the differential operator $\rho$ given by
\eqref{oper}. This result has been independently obtained by A.~Varchenko
by other methods. The appearance of Miura transformation in Gaudin models
has also been studied in \cite{Ush} from a different point of view.

2. If we allow singular points corresponding to other elements of the Weyl
group, we obtain more complicated equations of the type \eqref{beq} for
these points. These equations can also be interpreted as the equations on
the existence of a null-vector in the Wakimoto module corresponding to the
expansion of $\chi(t)$ at such a point. If the equations hold, then the
null-vector exists, and we can construct an eigenvector of the Gaudin
hamiltonians in the same way as in \cite{FFR}, \S~5.

Presumably, for generic $z_1,\ldots,z_N$ all eigenvectors of the Gaudin
hamiltonians correspond to the simplest degeneracies of the opers,
i.e. only those labeled by the simple reflections from the Weyl
group. But for special values of $z_1,\ldots,z_N$, some of the eigenvectors
may correspond to other degeneracies.

3. When $X$ is of genus one, the ${\mc D}$--modules obtained by geometric
Langlands correspondence can also be described explicitly. As special cases
one obtains the Calogero-Sutherland systems with elliptic potential
\cite{ER}.

\section{Drinfeld's construction and Sklyanin's separation of variables.}

V.~Drinfeld has given another construction of geometric Langlands
correspondence in the case of $GL_2$ \cite{Dr}. Recall from \secref{glob}
that the goal is to attach to a two-dimensional $l$--adic representation
$\sigma$ of $\pi_1(X)$, an automorphic function on $G(F)\backslash
G({\bolda})/K$. Such a function can be viewed as a $G(F)$--invariant
function on
$G({{\bolda}})/K$. H.~Jacquet and R.P.~Langlands \cite{JL} constructed a
function on $G({{\bolda}})/K$, which is $B(F)$--invariant, where $B$ is the
Borel subgroup of $G$. It remains then to establish that it is actually
$G(F)$--invariant. H.~Jacquet and R.P.~Langlands derived this from the
functional equation on the corresponding $L$--function (which they had
assumed from the beginning).

V.~Drinfeld's approach \cite {Dr} (cf. also \cite{Laumon}) was to associate
a perverse sheaf to the function constructed by Jacquet-Langlands, and show
that this sheaf is constant along the fibers of the projection
$G(F)\backslash G({{\bolda}})/K \arr B(F)\backslash G({{\bolda}})/K$.

\subsection{Whittaker function.} The Jacquet-Langlands function is
constructed using the {\em Whittaker model}. Let $G$ be a simple Lie
group. Choose a non-trivial additive character $\Phi: {{\bolda}}/F \arr
\Ql^\times$. Let $W$ be the representation of $G({{\bolda}})$ in the space of
locally constant functions $f: G({{\bolda}}) \arr \Ql$, such that
\begin{equation}    \label{whittaker}
f(u^{-1} x) = \sum_{i=1}^l \Phi(u_i) f(x)
\end{equation}
for any element $u$ of the unipotent subgroup $N({{\bolda}}) \simeq {\bolda}$.
 Here $u_i$ is the element of $u$ corresponding to the $i$th simple
root of $G$. A Whittaker model for a $G({{\bolda}})$--module $M$ is by
definition a non-trivial homomorphism $M \arr W$.

Recall that by local Langlands correspondence we can associate to a
homomorphism $\sigma: \pi_1(X) \arr G^L$ a collection of
unramified representations $\pi_x$ of the groups $G_x, x\in X$. It is known
that the $G({{\bolda}})$--module $\otimes_{x\in X} \pi_x$ has a unique
Whittaker model if all factors $\pi_x$ are non-degenerate \cite{GelK}. The
image of the vector $\otimes_{x\in X} v_x$ under the homomorphism to $W$ is
a function $\phi_\sigma$ on $G({{\bolda}})$, which is called {\em Whittaker
function}. This function is right invariant with respect to $K$ and
satisfies equation \eqref{whittaker}.

It is easy to see that $N({\mc K}_x) \backslash G_x/G({\mc O}_x)$ is
isomorphic to the root lattice of $G$ and hence to the weight lattice
$P^\vee$ of the Langlands dual group $G^L$. Therefore $N({{\bolda}})
\backslash G({{\bolda}})/K$ is isomorphic to the set $\on{Div}_{P^\vee}$
of $P^\vee$--valued divisors on $X$. Hence the Whittaker function
$\phi_\sigma$ can be considered as a function on $\on{Div}_{P^\vee}$. It
turns out \cite{Shalika} that its values on non-effective divisors are
equal to $0$, and
\begin{equation}    \label{shalika}
\phi_{\sigma}\left( \sum_i \la_i [x_i] \right) = \prod_i \on{tr}
\; \sigma(\on{Fr}_{x_i})|_{V_{\la_i}}, \quad \quad \la_i \in P^\vee_+,
\end{equation}
up to a power of $q$. Here $V_\la$ is the finite-dimensional
representation of $G^L$ of highest weight $\la \in P^\vee_+$.

\subsection{Geometric interpretation.}
\label{geomint} From now on we put $G=GL_2$. Then $\sigma$
can be considered as an
irreducible two-dimensional representation of $\pi_1(X)$. This
representation corresponds to a rank two locally constant sheaf
$E_\sigma$ on $X$, cf. e.g. \cite{Milne}. Since $P^\vee = \Z$ in this case,
$\on{Div}_{P^\vee}$ is just the set of divisors on $X$. The function
$\phi_\sigma$ is zero away from the subset $\on{Div}_+$ of effective
divisors, and its values on $\on{Div}_+$ are given by formula
\eqref{shalika}.

The set of effective divisors on $X$ of degree $m$ is the set of
points of the $m$th symmetric power of $X$, $S^m X$. Consider the
irreducible perverse sheaf $E_\sigma^{(m)} = \left( \pi_*
E_\sigma^{\boxtimes m} \right)^{S_m}$ on $S^m X$, where $\pi: X^m \arr
S^m X$ is the projection. The stalk of $E_\sigma^{(m)}$ over a divisor
$\sum_i n_i [x_i]$, where $x_i$'s are distinct, is $\otimes_i S^{n_i}
E_{\sigma,x}$. Let $\phi_\sigma^m$ be the restriction of $\phi_\sigma$
to $S^m X$. It follows from \eqref{shalika} that the value of
$\phi_\sigma^m$ at $d \in S^m X$ is given by the trace of the
Frobenius element of $d$ on the stalk of $E_\sigma^{(m)}$ over $d$.

Thus, geometrically, the passage from $\sigma$ to the Whittaker function is
the passage from $E_\sigma$ to $E^{(m)}_\sigma, m>0$. This is the first
step of Drinfeld's construction \cite{Dr}.

Now consider the function
\begin{equation}    \label{jl}
\wt{g}_\sigma(x) = \sum_{a\in F^\times} \phi \left(
\begin{pmatrix}
1 & a \\
0 & 1
\end{pmatrix} x \right).
\end{equation}
Formula \eqref{whittaker} implies that this function is left invariant
with respect to $B(F)$, and hence it gives rise to a function $g_\sigma$ on
$B(F)\backslash G({{\bolda}})/K$. This is the Jacquet-Langlands function.

Drinfeld has interpreted the space $B(F)\backslash G({{\bolda}})/K$ as
the set of points of the moduli space $\M_{2,1}$ of rank two bundles
on $X$ with a rank one subbundle. An element of $\M_{2,1}$ can be
considered after twisting by an appropriate line bundle as an
extension
\begin{equation}    \label{extension}
0 \larr {\mc O} \larr {\mc L}_2 \larr {\mc L}_1 \larr 0.
\end{equation}
Let $\M_{2,1}^n$ be the component of $\M_{1,2}$, which consists of
extensions \eqref{extension} where ${\mc L}_1$ and hence ${\mc L}_2$ are
of degree $n$. There is a map $j_n^\vee: \M_{2,1}^n \arr Jac_n$ to the
$n$th Jacobean variety $Jac_n$ of $X$, which sends an extension
\eqref{extension} to the class of ${\mc L}_1$. If $n > 2g-2$, $j_n^\vee$
is a projective fibration with the fiber $\Pro H^1(X,\on{Hom}({\mc
L}_1,{\mc O}))$.

The dual fibration is isomorphic to $S^m X$, where $m=n+2g-2$. Indeed,
there is natural map $j_n: S^m X \arr Jac_n$, which sends a divisor $d \in
S^m$ to ${\mc O}(d) \otimes \Omega^{-1} \in Jac_n$, where $\Omega$ is the
canonical bundle. The fiber of the map $j_n$ over ${\mc L}_1$ is $\Pro
H^0(X,{\mc L}_1 \otimes \Omega)$, which is dual to the fiber of
$j_n^\vee$ over ${\mc L}_1$, by the Serre duality.

For a pair of dual projective bundles over the same base one can define
{\em Radon transform}, which maps functions on one of them to functions on
the other. As shown in \cite{Dr}, formula \eqref{jl} actually means that
the restriction $g_\sigma^n$ of $g_\sigma$ to $\M_{2,1}^n$ is the Radon
transform of the function $\phi_\sigma^m$ on $S^m X$.

One can define the Radon transform on perverse sheaves,
cf. e.g. \cite{Br,KS,Laumon}. It follows from a theorem of P.~Deligne,
cf. \cite{Dr}, that the Radon transform of $E_\sigma^{(m)}$ is an
irreducible perverse sheaf on $\M_{2,1}^n$, which we denote by ${\mc
G}_\sigma^n$. The sheaf ${\mc G}_\sigma^n$ is the geometric object, which
replaces the function $g_\sigma^n$ in the sense that the value of
$g_\sigma^n$ at a point $p$ of $\M_{2,1}^n$ is given by the trace of the
Frobenius element of $p$ on the fiber of ${\mc G}_\sigma^n$ over $p$.

$$\begin{array}{ccccccc} \; & \; & \M_{2,1}^n & \; & \; & \; & S^m X \\ \;
& \stackrel{i_n}{\swarrow} & \; & \stackrel{j_n^\vee}{\searrow} &
\; & \stackrel{j_n}{\swarrow} \\ \M_2^n & \; & \; & \; & Jac_n & \; &
\;
\end{array}$$

\subsection{Descent.} There is a natural projection $i_n: \M_{2,1}^n \arr
\M_2^n$ to the moduli space of rank two bundles on $X$ of degree $n$, which
maps an extension \eqref{extension} to the class of ${\mc L}_2$. This is
of course just the geometric realization of the projection $G(F)\backslash
G({{\bolda}})/K \arr B(F)\backslash G({{\bolda}})/K$. The fiber of $i_n$
over ${\mc L}_2 \in \M_2^n$ is the set of non-vanishing sections of ${\mc
L}_2$ up to a non-zero multiple, which is the complement of a divisor in
$\Pro H^0(X,{\mc L}_2)$.

\smallskip
\noindent{\em Remark.} Here we consider $i_n$ as a map of stacks. But we
could also view it as a rational map from the moduli space of semi-stable
pairs (bundle,section) to the moduli space of semi-stable bundles. In this
setting, beautiful resolutions of this map have been given by A.~Bertram
\cite{Ber} and M.~Thaddeus \cite{Th} (cf. also \cite{FS}).\qed
\smallskip

Drinfeld proves that ${\mc G}_\sigma^n$ is constant along the generic
fiber of $i_n$. This implies that ${\mc G}_\sigma^n$ is a pull-back of a
perverse sheaf ${\mc F}_\sigma^n$ on $\M_2^n$ for all $n>2g-2$. After that
Drinfeld constructs the automorphic function as the function attached to
the sheaves ${\mc F}_\sigma^n$; using the Hecke operators this function
can be uniquely extended to the whole $\M_2$ \cite{Dr}.

\subsection{Comparison of two constructions.} Since the construction
outlined above is purely geometric, it makes sense over the field of
complex numbers as well, and it establishes geometric Langlands
correspondence for $GL_2$ in the sense of \secref{ho} (more precisely, in
one direction, while over $\F$ this has also been done in the opposite
direction, cf. \cite{Dr:shtuka}). It is interesting to compare it with the
construction outlined in \secref{ho}. Over the field of complex numbers,
one can switch from perverse sheaves to ${\mc D}$--modules using the
Riemann-Hilbert correspondence, cf. e.g. \cite{KS}.

Let us also switch from $GL_2$ to $SL_2$: we have $\M_{SL_2}(X) \simeq
\M_2^n$ for even $n$. Let us take as ${\mc F}_\sigma^n$ the ${\mc
D}$--module $\Delta_\rho$ obtained from a particular projective connection
$\rho$, whose monodromy representation coincides with $\sigma: \pi_1(X)
\arr PGL_2$. Now take its inverse image $i_n^* \delta(\rho)$ on
$\M_{2,1}^n$ and apply the inverse Radon transform. This gives us a ${\mc
D}$--module on $S^m X$, where $m=n+2g-2$. The anticipated equivalence of
the two constructions simply means that the latter is just the ${\mc
D}$--module corresponding to $E_\sigma^{(m)}$ by the Riemann-Hilbert
correspondence (i.e. that $E_\sigma^{(m)}$ is the sheaf of solutions of
this ${\mc D}$--module). Clearly, this would imply that $\Delta_\rho$ is
automorphic and depends on $\sigma$, but not on $\rho$, and that ${\mc
G}_\sigma$ is constant along the fibers.

In order to establish the equivalence of the two constructions, it is
important to understand what is the analogue of the Whittaker model over
the complex field. The existence of the Whittaker model over a finite field
is equivalent to the existence of a non-trivial coinvariant of the
unipotent subgroup $N$ in $M \otimes \Phi^*$, where $M$ is a $G$--module
and $\Phi^*$ is the character of $N$ dual to $\Phi$.

Now let $M$ be a $\G$--module from the category ${\mc O}^0_{{\rm crit}}$,
and $\chi$ be a character of $L\n_+$. A naive analogue of the local
Whittaker model is the space of coinvariants of $M \otimes \chi$ with
respect to the Lie subalgebra $L\n_+$. But this space is zero, because the
subalgebra $\n_+ \otimes \C[t]$ of $L\n_+$ acts locally nilpotently on
$M$. The correct analogue is the {\em semi-infinite} cohomology of $L\n_+$
with coefficients in $M \otimes \chi$ -- the quantum Drinfeld-Sokolov
reduction of $M$ introduced in \cite{FF:ctr}.

The analogue of the global Whittaker model should be a semi-infinite
localization functor assigning a ${\mc D}$--module on each $S^n X$ to
a $\G(\bolda)$--module $\otimes'_{x\in X} M_x$. In the case of $\su$
the fiber of this ${\mc D}$--module at the divisor $\sum_{x\in X} n_x
[x] \in \on{Div}_+$ should be the tensor product over all $x \in X$
of the semi-infinite cohomologies of $L\n_{+,x} = \n_+ \otimes
\C((t_x))$ with coefficients in $M_x \otimes \chi_x$, where $\chi_x$
is a character of $L\n_{+,x}$ vanishing at the origin up to order
$n_x$.

We can apply this functor to $\otimes_{x\in X} V^{\rho_x}$, where
$\rho$ is a projective connection on $X$. In this case one can show
that the fiber of the corresponding ${\mc D}$--module on $S^n X$ at
the divisor $\sum_{x\in X} n_x [x]$ is isomorphic to $\otimes_x
S^{n_x} (\C^2)$. Therefore we believe that this ${\mc D}$--module is
isomorphic to the $n$th symmetric power of the ${\mc D}$--module on
$X$ defined by $\rho$.

\subsection{Genus zero case.} Let us consider the genus zero case allowing
parabolic structures at the points $z_1,\ldots,z_N$ and $\infty$; Drinfeld
has generalized his construction in \cite{Dr1} to allow for parabolic
structures.

Put $n=N-1$; $Jac_n$ is just one point, and we have a pair of dual
projective spaces: $\Pro^\vee=\Pro \on{Ext}({\mc
O}(z_1+\ldots+z_N+\infty),{\mc O})$ and $\Pro = \Pro
H^0(\pone,\Omega(z_1+\ldots+z_N+\infty))$. The latter is isomorphic to
$S^{N-1} \pone$: if $X_i$'s are the residues of one forms from $\Pro$ and
$y_1,\ldots,y_{N-1}$ are natural coordinates on $S^{N-1} \pone$, the
isomorphism is given by
\begin{equation}    \label{transition}
\sum_{i=1}^N \frac{X_i}{t-z_i} dt = r
\frac
{\prod_{j=1}^{N-1}
(t-y_j)}{\prod_{i=1}^N (t-z_i)} dt,
\end{equation}
i.e. by the passage from a one-form to its zeroes.

There is a natural map $\Pro^\vee \arr \M^{(N)}$. Namely, an element of
$\Pro^\vee$ defines an extension $0 \arr {\mc O} \arr {\mc L}_2 \arr
{\mc O}(z_1+\ldots+z_N+\infty) \arr 0$. Then the canonical map ${\mc O}
\arr {\mc O}(z_1+\ldots+z_N+\infty)$ provides us with an extension ${\mc
L}'_2$ of ${\mc O}$ by ${\mc O}$ with flags at
$z_1,\ldots,z_N,\infty$. The fiber of the map $\Pro^\vee \arr \M^{(N)}$ is
isomorphic to $\pone$ without $N+1$ points -- this is the set of sections
of ${\mc L}'_2$, which miss all flags, up to a scalar. There is a nice
geometric approximation of this map, due to M.~Thaddeus.

The ${\mc D}$--module $\Delta(\la_i,\mu_i)$ describes the Gaudin
equations $H_i \Psi(X_1,\ldots,X_N) =$
\newline
$\mu_i \Psi(X_1,\ldots,X_N),
i=1,\ldots,N$, which are equivalent to the equations
\begin{equation}    \label{initial}
S(t) \Psi = \left( \sum_{i=1}^N \frac{\la_i(\la_i+2)/4}{(t-z_i)^2} +
\sum_{i=1}^N \frac{\mu_i}{t-z_i} \right) \Psi, \quad t \in \C \backslash \{
z_1,\ldots,z_N \},
\end{equation}
where $S(t)$ is given by formula \eqref{su} and the action of the $i$th
copy of $\sw_2$ is given by formula \eqref{sl2}. The pull-back of this
system to $\Pro^\vee$ is the same system in which $x_i$'s should be
considered as the homogeneous coordinates on $\Pro^\vee$ dual to the
coordinates $X_i$ on $\Pro$.

The Radon transform is equivalent to the formal Fourier transform,
i.e. substituting $x_i \arr -\pa/\pa X_i$, and $\pa/\pa x_i \arr X_i$ in
formula \eqref{sl2}, so that the action of the $i$th copy of $\sw_2$
becomes
\begin{equation}    \label{newsl2}
e^{(i)} = -X_i^2 \frac{\pa}{\pa X_i} - (\la_i+2) X_i, \quad h^{(i)} = -2
X_i \frac{\pa}{\pa X_i} - (\la_i+2), \quad f^{(i)} = X_i.
\end{equation}
We now have to rewrite the system \eqref{initial} in the new coordinates
$y_j$'s using formula \eqref{transition}.

\subsection{Separation of variables.}
This has been done by E.~Sklyanin, using a clever trick
\cite{Skl:sep}. Observe that $S(t) = f(t) e(t) + \frac{1}{4} h(t)^2 -
\frac{1}{2} \pa_t h(t)$, where we use notation \eqref{ft}. After the
Fourier transform, $f(t) = \sum_{i=1}^N X_i/(t-z_i)$, and according to
formula \eqref{transition}, $f(t)|_{t=y_j} = 0$.

On the other hand, we have by \eqref{newsl2}: $$\frac{1}{2} h(t) = -
\sum_{i=1}^N \left( \frac{X_i \pa/\pa X_i}{t-z_i} - \frac{\la_i+2}{t-z_i}
\right).$$ But we find from \eqref{transition}:
$$\frac{\pa}{\pa y_j} = \sum_{i=1}^N \frac{X_i \pa/\pa X_i}{y_j-z_i}.$$
Hence
$$\frac{1}{2} h(t)|_{t=y_j} = -\nabla_{y_j} \equiv - \frac{\pa}{\pa
y_j} + \frac{1}{2} \sum_{i=1}^N \frac{\la_i+2}{y_j-z_i}.$$ Substituting
$t=y_j$ from the left into $S(t)$ and using formulas above, we obtain
\cite{Skl:sep}: $S(t)|_{t=y_j} = \nabla_{y_j}^2$.

Thus, substituting $t=y_j$ from the left into \eqref{initial} we obtain
the separated equations:
\begin{equation}    \label{separate}
\left( \nabla_{y_j}^2 - \sum_{i=1}^N \frac{\mu_i}{y_j-z_i} - \sum_{i=1}^N
\frac{\la_i(\la_i+2)/4}{(y_j-z_i)^2} \right) \Psi = 0, \quad
j=1,\ldots,N-1,
\end{equation}
on symmetric polynomials $\Psi(y_1,\ldots,y_{N-1})$. By a standard
interpolation argument, we see that the systems \eqref{initial} and
\eqref{separate} are equivalent.

Note that each of the equations \eqref{separate} has the form $(\nabla_t^2
- q(t)) \Psi=0$, where $\pa_t^2 - q(t)$ is the projective connection
\eqref{qt} from which the original ${\mc D}$--module $\Delta(\la_i,\mu_i)$
was obtained. The replacement of $\pa_t$ by $\nabla_t$ just means twisting
by a rank one local system, which does not change our $PGL_2$--oper.

The perverse sheaf, which corresponds to the ${\mc D}$--module on $S^{N-1}
\pone$ defined by the system \eqref{separate}, is isomorphic to
$E^{(N-1)}$, where $E$ is the sheaf of local solutions of the projective
connection $\pa_t^2 - q(t)$. Thus, the equivalence of Drinfeld's first and
Beilinson-Drinfeld second constructions amounts in genus zero to separation
of variables in the Gaudin equations. We will discuss this construction and
its generalizations in detail elsewhere.

The separation of variables has been found by Sklyanin in an attempt to
find an alternative to the conventional Bethe ansatz
\cite{Sklyanin,Skl:sep,Skl:quant}. Sklyanin's approach goes back to the
classical Hamilton-Jacobi theory of separation of variables in hamiltonian
systems. In the theory of non-linear equations, which are integrated by the
inverse scattering method, the separating variables are usually the zeroes
of the Baker function. Sklyanin's idea was that the function $f(t)$ above
should be viewed as an analogue of the Baker function, and hence the
equations should separate with respect to its zeroes. One of the advantages
of this approach is that one does not need the ``pseudo-vacuum'' $|0\ri$
for constructing eigenvectors, so representations without highest weight
can also be considered. It also gives the most straightforward answer to
the completeness problem.

Indeed, the Verma module $M_{\chi_i}$ over $\sw_2$ can be realized in the
space of polynomials $\C[X_i]$ with the action given by formula
\eqref{newsl2} with $\chi_i=-\la_i-2$. A solution of the system
\eqref{separate} is a product of solutions with respect to each of the
variables: $\Psi(y_1,\ldots,y_{N-1}) = \prod_{i=1}^{N-1} \psi(y_i)$. Thus,
a polynomial solution of \eqref{separate} should be of the form
$\Psi(y_1,\ldots,y_{N-1}) = \prod_{i=1}^{N-1} \prod_{j=1}^m
(y_i-w_j)$. Using formula \eqref{transition} we can pass to the
$X_i$--coordinates, and obtain an eigenvector of the Gaudin hamiltonians in
$\otimes_{i=1}^N M_{\chi_i}$. It coincides with the vector
$|w_1,\ldots,w_m\ri$ given by \eqref{bethevector} up to a scalar
multiple. Since the Gaudin equations on $\otimes_{i=1}^N M_{\chi_i}$ are
equivalent to equations \eqref{separate}, we obtain completeness of Bethe
ansatz for $\otimes_{i=1}^N M_{\chi_i}$. For other applications of the
separation of variables, cf. \cite{KH}.

There should exist a similar separation of variables for the generalized
(spin) Calo\-gero-Sutherland models defined on the moduli space of rank two
bundles over the torus with parabolic structures.

\subsection{Generalization to higher rank.}
G.~Laumon \cite{Laumon} suggested a geometric construction of Langlands
correspondence in the case of $GL_n$ generalizing V.~Drinfeld's first
construction for $GL_2$ and P.~Deligne's construction for $GL_1$. The idea
is to connect the moduli space of rank $r>1$ bundles with a section and the
moduli space of rank $r+1$ bundles with a section, by a pair of dual
projective fibrations of the kind given by the diagram above. This way one
continues the diagram further to the left and ultimately reaches the moduli
space of rank $n$ bundles. In order to construct a sheaf on the latter, one
has to apply a chain of Radon transforms to the sheaf $E^{(m)}$ on $S^m X$,
where $E$ is now a rank $n$ local system on $X$. The difficulty is to prove
analogues of the Deligne vanishing theorem (cf. \secref{geomint}), which
are necessary to insure that the Radon transforms do not spoil
irreducibility and perversity of sheaves, cf. \cite{Laumon}.

It would be interesting to compare Laumon's construction, more precisely,
its generalization to the bundles with parabolic structures, to the
Beilinson-Drinfeld construction in genus zero. This should provide a
separation of variables for Gaudin models of higher rank. For instance, in
the case of $SL_3$, the Gaudin system is a system of differential equations
in $3N$ variables, corresponding to an $\sw_3$--oper $\rho$, i.e. a third
order differential operator \eqref{oper}. We first have to apply a formal
Fourier transform with respect to $2N$ of them, and make a change of
variables analogous to \eqref{transition}. Then we have to apply another
Fourier transform with respect to $N$ variables, and apply another change
of variables. This should give us a system of identical equations of the
form $\rho \cdot \Psi = 0$. On the other hand, for the classical (and
partially for quantum) Gaudin models associated to $SL_n$, a different
scheme of separation of variables has been suggested in
\cite{Harnad,class}, and in \cite{Sklyanin} a separation of variables has
been given for the Calogero-Sutherland model corresponding to $SL_3$
and genus one. So there may exist another generalization of Drinfeld's
first construction to groups of higher rank.

\subsection{Quantization.}
As we already mentioned in the introduction, Sklyanin has found a quantum
deformation of the separation of variables of the Gaudin equations
\cite{Skl:quant}, in which the role of the differential equations of second
order \eqref{separate} is played by $q$--difference equations of second
order. This suggests that elements of the spectrum of the center $Z_q(\G)$
of a completion of the quantum affine algebra $U_q(\G)$ at the critical
level should be viewed as $q$--difference operators.

Explicit computation of the $q$--deformation of the Miura transformation
made in \cite{FR} shows that it is indeed so. For example, $Z_q(\su)$ is
generated by the Fourier components of a generating function $\ell(z)$,
which is a $q$--deformation of the Sugawara operator $S(z)$. The
$q$--deformation of the Miura transformation \eqref{mt} reads:
\begin{equation}    \label{map}
\ell(z) \longrightarrow \Lambda(qz) + \Lambda(zq^{-1})^{-1},
\end{equation}
where $\Lambda(z)$ is a $q$--deformation of $\chi(z)$. According to this
formula, $\ell(z)$ is a ``quantum trace'' of the conjugacy class of
$\on{diag}[\Lambda(z),\Lambda(z)^{-1}]$. Therefore elements of the
spectrum of $Z_q(\su)$ look very similar to the local Langlands parameters
described in \secref{ll}.

Formula \eqref{map} can be rewritten as follows: ${\mc D}_q^2 - \ell(z)
{\mc D}_q + 1 = ({\mc D}_q - \La(zq))({\mc D}_q - \La(zq)^{-1})$, where
$[{\mc D}_q \cdot f](z) = f(zq^{-2})$. This is a $q$--analogue of
\eqref{mt}. If $Q(z)$ is a solution of the $q$--difference equation $Q(zq)
= \La(z) Q(zq^{-1})$, then $({\mc D}_q + {\mc D}_q^{-1} - \ell(z)) Q(z) =
0$. This formula, which is a $q$--analogue of \eqref{separate}, coincides
with the equation which appears in quantum separation of variables
\cite{Skl:quant}. It plays a prominent role in statistical mechanics,
cf. \cite{Baxter}.

Quantum affine algebras may shed new light on the isomorphism between the
center at the critical level and a $\W$--algebra. Indeed, the $R$--matrix
of $U_q(\G)$ corresponding to a finite-dimensional representation coincides
with the $R$--matrix of the quantum Toda system associated to $\G$
\cite{JB}. The transfer-matrices of the Toda system are the integrals of
motion of the latter, and are closely related to the central elements of
$U_q(\G)$ at the critical level constructed in \cite{RS}, cf. \cite{FR}. On
the other hand, the classical $\W$--algebra $\W(\g^L)$ is just the algebra
of integrals of motion of the classical Toda system associated to $\g^L$,
and $\W(\g^L)$ contains integrals of motion of the Toda system associated
to $\widehat{\g^L}$, cf. \cite{FF:law,F:talk}. However, the appearance of
the Langlands dual Lie algebra here still remains a mystery.

\end{document}